\documentclass{CSML}

\def\dOi{12(2:2)2016}
\lmcsheading%
{\dOi}
{1--25}
{}
{}
{Jun.~15, 2015}
{Apr.~19, 2016}
{}

\ACMCCS{[{\bf Theory of computation}]: Models of
computation---Computability---Turing
machines---Abstract machines; Computational complexity and cryptography;
Logic---Logic and verification==-Constructive mathematics.}

\usepackage{amssymb, amsmath, amsthm, fancyhdr, fancybox, hyperref}

\begin{document}
\renewcommand{\phi}{\varphi}

\newcommand{\qeD}{{\nobreak\hfil\penalty50
  \hskip1em\hbox{}\nobreak\hfil$\sqr{8}{8}$
  \parfillskip=0pt \finalhyphendemerits=0 \par}
  \bigskip}
\newcommand{\prf}{\noindent{\sc Proof.}\enspace}
\newcommand{\Prf}[1]{\bigbreak\noindent{\sl Proof #1.}\enspace}
\newcommand{\myprime}{{\rm prime }}
\newcommand{\chain}{{\rm chain }}
\newcommand{\unchain}{{\rm unchain }}
\newcommand{\Iff}{\Leftrightarrow}

\newtheorem{prp}{Proposition}
\newtheorem{lemm}{Lemma}
\newtheorem{subclm}{Subclaim}
\newtheorem{assume}{Assumption}
\newtheorem{rmrk}{Remark}

\title[Non-Obfuscated Yet Unprovable Programs]
      {Non-Obfuscated Unprovable Programs \& Many Resultant Subtleties\rsuper*}

\author[J.~Case]{John Case}
\address{Computer \& Information Sciences,
	   University of Delaware,
	  Newark, DE 19716 USA}
\email{\{case, mralston\}@udel.edu}

\author{Michael Ralston}
\address{\vspace{-18 pt}}

\begin{abstract}

The \emph{International Obfuscated C Code Contest}
was a programming contest for the most creatively obfuscated yet
succinct C code.
By \emph{contrast}, an interest herein is in programs which are,
\emph{in a sense}, \emph{easily} seen to be correct,
but which can\emph{not} be proved correct in pre-assigned, computably
axiomatized, powerful, true theories {\bf T}.  A point made by our
first theorem, then, is that, then,
\emph{un}verifiable programs need \emph{not} be obfuscated!

The first theorem and its proof is followed by a motivated, concrete
example based on a remark of Hilary Putnam.

The first theorem has some non-constructivity in its statement and
proof,
and the second theorem implies some of the non-constructivity is
inherent.
That result, then, brings up the question
of whether there is an acceptable programming system (numbering) for
which
some non-constructivity of the first theorem disappears. The third
theorem
shows this is the case, but for a subtle reason explained in the text.
This
latter theorem has a number of corollaries, regarding its 
acceptable
programming system, and providing some surprises and
subtleties about proving its program properties (including
universality, and
the presence of the composition control structure).
The next two theorems provide acceptable systems with contrasting
surprises
regarding proving universality in them.  Finally the next and last
theorem
(the most difficult to prove in the paper) provides an acceptable
system with
some positive and negative
surprises regarding verification of its true program properties: the existence
of the control structure composition is provable for it, but anything
about
true I/O-program equivalence for syntactically
unequal programs is not provable.

\end{abstract}


\titlecomment{{\lsuper*}Some of the results in the present paper were presented
at
the \emph{Twelfth Asian Logic Conference}, Wellington, NZ, 2011, under
the
title, Non-Obfuscated Yet Unprovable Programs.}

\renewcommand{\phi}{\varphi}

\newcommand{\formula}[1]{{\ll #1 \gg}}
\newcommand{\fqE}{{ \ll {\sf E} \gg }}
\newcommand{\fqEp}{{ \ll {\sf E}' \gg}}
\newcommand{\PA}{{\bf PA}}
\newcommand{\compose}{\mathclose{\circ}\mathclose}
\newcommand{\sO}{{\mathcal O}}
\newcommand{\ifthenelseforus}{\mbox{\rm if-then-else}}
\newcommand{\krt}{\mbox{\rm krt}}
\newcommand{\Implies}{\Rightarrow}
\newcommand{\pad}{\mbox{\rm pad}}
\newcommand{\KRT}{\mbox{\bf KRT}}
\newcommand{\Lineartime}{{\bf LinearTime}}
\newcommand{\T}{\mbox{{\bf T}}}
\newcommand{\proves}{\vdash}
\newcommand{\notproves}{\mathrel{\not\,\vdash}}
\newcommand{\AND}{\ \wedge \ }
\newcommand{\OR}{\ \vee \ }
\newcommand{\N}{\mathbb{N}}
\newcommand{\natnum}{\mathbb{N}}
\newcommand{\sqdot}{\rule{0.5mm}{0.5mm}}
\newcommand{\bigsqdot}{\rule{1.5mm}{1.5mm}}
\newcommand{\lam}[1]{\lambda #1 \,\sqdot\,}
\newcommand{\converges}{\mathclose{\hbox{$\downarrow$}}}
\newcommand{\diverges}{\mathclose{\hbox{$\uparrow$}}}
\newcommand{\set}[1]{\{#1\}}
\newcommand{\pair}[1]{\langle #1 \rangle}
\newcommand{\Kbar}{\overline{K}}
\newcommand{\symmdiff}{\bigtriangleup}
\newcommand{\Quads}[1]{\hspace*{#1em}}
\newcommand{\Or}{\ \vee \ }
\newcommand{\existsio}{\mathord{\stackrel{\infty}{\exists}}}

\newcommand{\sqr}[2]{{\vcenter{\hrule height.#2pt
	\hbox{\vrule width.#2pt height#1pt \kern#1pt
	   \vrule width.#2pt}
	\hrule height.#2pt}}}

\renewcommand{\Qed}[1]{{\nobreak\hfil\penalty50
  \hskip1em\hbox{}\nobreak\hfil$\sqr{8}{8}$\enspace\sc #1
  \parfillskip=0pt \finalhyphendemerits=0 \par}
  \bigskip}

\newtheorem{defn}{Definition}
\newtheorem{thrm}{Theorem}

\newcommand{\cardb}[1]{\mid\mathclose{#1}\mathclose{\mid}}
\newcommand{\floor}[1]{\lfloor#1\rfloor}
\newcommand{\ceil}[1]{\lceil#1\rceil}
\newcommand{\eqdef}{\stackrel{\rm def}{=}}
\newcommand{\oldphitm}{\varphi^{\rm TM}}
\newcommand{\phitm}{\varphi}
\newcommand{\Phitm}{\Phi^{\rm TM}}
\newcommand{\Ptimek}{{\mathcal P}\mbox{\rm time}_k}
\newcommand{\Ptime}{{\mathcal P}\mbox{\rm time}}
\newcommand{\Pspace}{{\mathcal P}\mbox{\rm space}}
\newcommand{\Ltime}{{\mathcal L}\mbox{\rm time}}
\newcommand{\comp}{\mbox{\rm comp}}
\newcommand{\fncmp}{\mbox{\rm \footnotesize cmp}}
\newcommand{\cmp}{\mbox{\rm cmp}}
\newcommand{\fncomp}{\mbox{\rm \footnotesize comp}}
\newcommand{\trans}{\mbox{\rm trans}}
\newcommand{\fntrans}{\mbox{\rm \footnotesize trans}}
\newcommand{\bound}{\mbox{\rm bound}}
\newcommand{\Phitmsd}{\Phi^{\rm SlowedDownTM}}
\newcommand{\pads}{{\mbox{{\rm\scriptsize pad}}}}

\maketitle

\section*{Introduction}\label{intro}

The \emph{International Obfuscated C Code Contest} (see the Wikipedia
entry)
was a
programming contest for the most creatively obfuscated C code, held
annually
between 1984 and 1996, and thereafter in 1998, 2000, 2001, 2004, and
2006.

In many cases, the winning programmer did something simple in such an
obscure but succinct way that it was hard for other (human)
programmers to see how his/her code actually worked.

By \emph{contrast}, our first interest herein is in programs which are,
\emph{in a sense}, \emph{easily} seen to be correct,
but which can\emph{not} be proved correct in pre-assigned, computably
axiomatized, powerful, true theories {\bf T}.  A point is that, then,
\emph{un}verifiable programs need \emph{not} be obfuscated!

Our first theorem (Theorem~\ref{thm:exists-g-w-domain-two} in
Section~\ref{main-result} below)
entails:
for \emph{any} deterministic, multi-tape Turing Machine (TM) program
$p$,
there will
be an \emph{easily seen equivalent} such TM program $q$ \emph{almost
(i.e.,
within small, linear factors) as fast and succinct as $p$}, but this
equivalence will \emph{not} be provable in {\bf T}.

A point of the
just mentioned, small, linear factors is that the \emph{un}provability
is \emph{not} based on some huge (or at least non-linear) growth in
run-time and/or program size in passing from $p$ to $q$.  In fact we'll
see in the proof of the first theorem that $q$ will be like $p$ except
that
$q$, in effect, encapsulates $p$
in a top-level if-then-else with: 1.~$p$ being the
else-part and 2.~the succinct, linear-time testable
if-condition being easily seen never to come true
 (but with this never coming true being \emph{un}provable in $\T$).

A motivated, concrete, special case, based on a remark in
Putnam \cite{Put:i:80}, will be presented (also in
Section~\ref{main-result} below).

As will be seen,
the first theorem and its proof have some
non-constructivity, and, with Theorem~\ref{thm:no-g-for-phi} in
Section~\ref{constructivity} below, some of this
non-constructivity is seen to be inherent.

Considered next, in Section~\ref{subtleties},
is whether the just mentioned non-constructivity goes away for some
acceptable
programmming systems (numberings).  The answer
(Theorem~\ref{thm:exists-a-g})
is affirmative, but for pleasantly \emph{subtle}
reasons spelled out in the section.
This
latter theorem has a number of corollaries
(Corollaries~\ref{cor:t-not-prove-psi-eq-phi},~\ref{cor:not-krt-and-u},
and~\ref{cor:t-not-prove-eps})
regarding its acceptable
programming system, and they provide some surprises and
subtleties about proving its program properties, including in
Corollary~\ref{cor:not-krt-and-u}, about universality and
the presence of the composition control structure.

Section~\ref{subtleties} makes up most of the paper.
Also within it
(in Sections~\ref{univ} and~\ref{comp})
are presented a number of positive and negative
surprises regarding verification of true program properties.

In Section~\ref{univ} 
Theorems~\ref{thm:t-prove-one-u}
and~\ref{thm:t-not-prove-u-in-theta} provide respective acceptable
systems with contrasting surprises
regarding proving universality in them.
Of course any acceptable system has infinitely many universal programs,
but Theorem~\ref{thm:t-prove-one-u} provides an acceptable system in
which exactly one of these universal programs is provably so.
By contrast, Theorem~\ref{thm:t-not-prove-u-in-theta} provides a
different
acceptable system with \emph{no} program which is \emph{provably}
universal.

Finally, in Section~\ref{comp}, the next and last theorem
(Theorem~\ref{thm:t-prove-composition-not-equality}), which is
the most difficult to prove in the paper, also
provides an acceptable system with some positive and negative
surprises regarding verification of its true program properties: for this
acceptable system, the existence
of the control structure composition is provable for it, but anything
about
true I/O-program equivalence for syntactically
unequal programs is not provable.


   \section{Mathematical Preliminaries}\label{prelims}

      \subsection{Complexity-Bounded Computability}
		    \label{complexity-bounded}

Let $\oldphitm$ be the \emph{efficiently} laid out and
G\"odel-numbered acceptable programming system (numbering) from
\cite[Chapter~3 \&~Errata]{Roy-Cas:b:94} and which is based on
deterministic
\emph{multi-tape Turing Machines} (with base~two I/O).\footnote{In
general,
the {\em acceptable\/} programming systems
\cite{Rog:j:58,Rog:b:67,Mac-You:b:78,Ric:th:80,Ric:j:81,Roy:b:87}
can be characterized as those programming
systems for all the $1$-argument partial computable functions:
$\N \rightarrow \N$  which are
inter-compilable with the natural system
$\oldphitm$. Rogers \cite{Rog:j:58,Rog:b:67} characterized
the acceptable systems as those with universality and for which
Kleene's
S-m-n
holds. This latter is more than enough to get (not necessarily
efficient)
recursion theorems in acceptable systems.}
Its programs are
named by all the numbers in $\natnum = \set{0,1,2,\ldots}$.
$\oldphitm_p$ is the partial computable function $\N \rightarrow \N$
computed by $\oldphitm$-program (number) $p$.
The numerical naming just mentioned does
\emph{not} feature prime powers and factorization, but, instead, is a
\emph{linear-time computable and invertible coding}.
Let $\Phitm$ be the
corresponding step-counting Blum Complexity Measure
\cite{Blu:j:67:complexity}.
$(\oldphitm, \Phitm)$ is a \emph{base} model for deterministic run time
costs.  $\oldphitm$'s superscript is  awkward when $\oldphitm$ 
is employed in subscripts, so, from this point on, we will write
$\oldphitm$ as simply $\phi$.

Herein, we will use the linear-time computable and invertible
pairing function $\pair{\cdot,\cdot}$ from \cite{Roy-Cas:b:94}:
the binary representation of $\pair{x, y}$ is (by definition)
an interleaving of the binary
representations of $x$ and $y$ where we alternate $x$'s and $y$'s digits
and start on
the right with the least most significant $y$ digit. For example, 
$\pair{15,2} = 94$ ---
since $15 = 1111$ (binary), $2 = 0010$ (binary), 
and $94 = 10101110$ (binary).
This function, clearly then,
maps all the \emph{pairs} of elements of $\N$ 1-1, onto $\N$.
We also employ this notation, based on
iterating, $\pair{\cdot,\cdot}$, as in \cite{Roy-Cas:b:94},
to code also triples, quadruples, \ldots
of elements of $\N$ 1-1, onto $\N$: for all $n>2$, and all
$x_1, \ldots, x_n+1, \pair{x_1, \ldots, x_n+1} =
\pair{x_1, \pair{x_2, \ldots, x_n+1}}$.  
These functions also clearly satisfy the following
\begin{lem}\label{lem-pairing} {\ }
\begin{enumerate}
\item $\pair{x_1, \ldots, x_n}$ is odd implies $x_n$ is odd;
\item $\lam{x_1, \ldots, x_n}\pair{x_1, \ldots, x_n}$ is monotonically
increasing in each of its arguments; and,
\item for all $x_1, \ldots, x_n, \max(x_1, \ldots, x_n)  \leq 
\pair{x_1, \ldots, x_n}$.
\end{enumerate}
\end{lem}
For example,
in the proof of Theorem~\ref{thm:t-prove-composition-not-equality} below,
the just above lemma will see explicit and implicit application.

$\Lineartime$ is the class of functions:
$\N \rightarrow \N$ each computable by some $\phitm$-program running
within
a $\Phitm$-time bound linear in the \emph{length} of its
base-two expressed argument.  Of course by means of the iterated
$\pair{\cdot,\cdot}$
function defined just above, we can and sometimes will speak of
multi-argument functions as being (or not being) in $\Lineartime$.

For $k\in\N$, $k$ could be a
numerically named program of $\phitm$ or just a datum. We let
$\cardb{k} = $ the \emph{length} of $k$, where $k$ is written \emph{in
binary}.
We can write this length as $(\ceil{\log_2 (k+1)})_+$, \emph{where}
$(\cdot)_+$ turns $0$ into $1$; else, leaves unchanged.\footnote{This
formula can be derived as the minimum number of whole bits needed
to store any one of the $k+1$ things $0$ through
$k$, \emph{except} that the case of $k=0$ needs only $0$ bits; however,
a single $0$ has length $1$.

This and
more general use of $(\cdot)_+$ \emph{also} helps to deal with the fact
that
zero values can cause trouble for $\sO$-notation ($\sO$-notation is
explained in \cite{CLRS:2001}).
A problem comes with complexity bounds
of more than one argument.      Jim Royer gave the following example of
two
functions mapping pairs from $\N$,
 $f(m,n) = (m \cdot n)$ \& $g(m,n) = (m+1)\cdot(n+1)$.  Suppose, as
\emph{might} be expected,
$g$ is $\sO(f)$.  Then there are positive $a,b$ such that,
for each $m,n \in \N$, $g(m,n) \leq
a \cdot f(m,n) + b $.  Then we have, for \emph{each}
$n$, $n+1 \leq a \cdot f(0,n) + b = b$, a contradiction.
\emph{However}, $g$ \emph{is} $\sO((f)_+)$.}

Rogers \cite{Rog:b:67} uses the terms `converges' for computations
which
halt and provide output and `diverges' for those that do not.
Herein we use the respective notations (due to Albert Meyer)
$\converges$ and $\diverges$ in place of those terms of Rogers.

From \cite[Lemma~3.14]{Roy-Cas:b:94},
there are small positive $a\in\natnum$ and function
$\ifthenelseforus \in \Lineartime$
such that, for all $p_0,p_1,p_2,x \in \natnum$,
   \begin{equation}
\phitm_{\ifthenelseforus(p_0,p_1,p_2)}(x) = 
   \begin{cases}
\phitm_{p_1}(x), & \mbox{if }\phitm_{p_0}(x)\converges \neq 0;\cr
\phitm_{p_2}(x), & \mbox{if }\phitm_{p_0}(x)\converges = 0;\cr
      \diverges, & \mbox{otherwise};\cr
   \end{cases}
   \end{equation}
and
      \begin{equation}\label{condefficienttoo}
 \Phitm_{\ifthenelseforus(p_0,p_1,p_2)}(x) \leq 
   \begin{cases}
a\cdot(\Phitm_{p_0}(x)+\Phitm_{p_1}(x))_+, &\mbox{if }
\phitm_{p_0}(x)\converges \neq 0;\cr
a\cdot(\Phitm_{p_0}(x)+\Phitm_{p_2}(x))_+, &\mbox{if }
\phitm_{p_0}(x)\converges = 0;\cr
      \diverges, & \mbox{otherwise}.\cr
\end{cases}
   \end{equation}\medskip

\noindent Essentially
from (the $k,m=1$ case of)
\cite[Theorem~4.8]{Roy-Cas:b:94}, we have the following constructive,
\emph{efficient},
and {\em parametrized\/} version of \emph{Kleene's} 2nd (not Rogers')
Recursion Theorem \cite[Page~214]{Rog:b:67}. \\
There are small positive $b\in\natnum$ and function
$\krt \in \Lineartime$ such that,
for all \emph{parameter values $p$, tasks $r$, inputs $x$} $\in
\natnum$:
   \begin{equation}\label{krt}
\phitm_{\krt(p,r)}(x) =
\phitm_{r}(\pair{\krt(p,r), p, x});
   \end{equation}
and
      \begin{equation}\label{complexity-task}
\Phitm_{\krt(p,r)}(x) \leq 
b\cdot(\cardb{p} + \cardb{r} + \cardb{x} +
\Phitm_{r}(\pair{\krt(p,r),p,x})).
   \end{equation}
    Intuitively, above in~(\ref{krt}),
on the left-hand side, the $\phitm$-program
\emph{$\krt(p,r)$ has $p,r$ stored inside}, and, on $x$, it:
\emph{makes a self-copy} (in linear-time),
\emph{forms $y=\pair{\mbox{self-copy},p,x}$} (in linear-time), and
\emph{runs task $r$ on this $y$}.  From~(\ref{complexity-task}) just
above, for each $p,r$,
any super-linear cost of running $\phitm$-program $\krt(p,r)$
on its input is from the running of $\phitm$-task $r$ on its
linear-time
producible input.

      \subsection{Computably Axiomatized, Powerful, True Theories}
			       \label{applied-logics}

Let ${\bf T}$ be a \emph{computably axiomatized}
first order (fo) theory extending fo Peano arithmetic
$(\PA)$ \cite{Men:b:09,Rog:b:67} --- \emph{but with numerals
represented in base~two to avoid size blow up from unary
representation}
(see \cite[Page~29]{Bus:th:86})\footnote{%
Lets
suppose $\overline{0}$ is $\PA$'s numeral for zero and that $S$ is
$\PA$'s
symbol for the successor function on $\N$.  In effect,
in, e.g., \cite{Men:b:09}, the
numeral $\overline{n}$
for natural number $n$ is $S^{(n)}(\overline{0})$, where
$S^{(0)}(\overline{0})=\overline{0}$ and 
$S^{(n+1)}(\overline{0})=S(S^{(n)}(\overline{0}))$ --- 
featuring iterated composition of $S$s.  This is a base~{\em one\/}
representation.
Note that the length of
{\em this\/} $\overline{n}$ is $\sO(n)$ which is
$\sO$ of $2^{\mbox{the symbol length of } \overline{n}}$ --- too high
for
feasible complexity.  However, the symbol length
for the binary representation of $n$ grows only linearly with $n$ ---
feasibly.

Based on \cite[Page~29]{Bus:th:86}, herein, by constrast with the just
above, we can define our numeral $\overline{n}$ for $n \in \N$ thus.
We suppose $\cdot$ is the symbol
for $\PA$'s multiplication over $\N$.
We let: $\overline{2} = S \compose S(\overline{0})$;
for $(n>1)$, $n \in \N$,
\begin{equation}
\overline{2n} = (\overline{2} \cdot \overline{n});
\end{equation}
and, for $n \in \N$,
\begin{equation}
\overline{(2n) + 1} = S(\overline{2n}).
\end{equation}
\emph{Then}, the length of $\overline{n}$ is in
$\sO$ of the symbol length of $\overline{n}$ --- feasible.
}
--- \emph{and which does not prove
(standard model for $\PA$ \cite{Men:b:09})
falsehoods expressible in f.o.~arithmetic}.

{\bf T} could be, \emph{for example}: fo Peano arithmetic $(\PA)$
itself,
the two-sorted fo Peano arithmetic permitting quantifiers over numbers
and
sets of numbers \cite{Rog:b:67,Sim:b:99}
(a second order arithmetic), Zermelo Frankel Set Theory with Choice
({\bf ZFC}) \cite{Hal:b:74:naive},
{\bf ZFC} $+$ ones favorite large cardinal
axiom
\cite{Roi:b:11,jech,Dra:b:74:large-cardinals,Kan:b:08:large-cardinals},
etc.

If ${\sf E}$ is an expression such as `the partial function computed by
$\phitm$-program number $p$ is total' and which is expressible in $\PA$ and
where $p$ is a particular element of $\N$, we shall write $\fqE$ to
denote a
typically {\em naturally\/} corresponding, fixed
standard cwff (closed well-formed formula) of $\T$
which (semantically) expresses ${\sf E}$ --- and where $p$ is expressed
as
the corresponding numeral in base~two (as indicated above).
We have that
\begin{quotation}
     \noindent
     \emph{if ${\sf E'}$ is obtained from ${\sf E}$ by substituting a
     numerical
     value $k$, then $\fqEp$ can be algorithmically
     obtained from $\fqE$ in linear-time in $(\cardb{\fqE} +
     \cardb{k})$.}
\end{quotation}\medskip

\noindent By \cite[Theorem~3.6 \& Corollary~3.7]{Roy-Cas:b:94} and their proofs,
the running of a {\em carefully crafted,
time-bounded, $\phitm$-universal simulation\/} up through time $t$
takes time a
little worse than exponential in $|t|$.
Early complexity theory, e.g., \cite{BCH:1969,Machtey:1973,Ladner75},
provided delaying tricks to achieve polynomial time.
From \cite[Theorem~3.20]{Roy-Cas:b:94} and its proof, the above
mentioned carefully crafted,
time-bounded universal simulation of any $\phitm$-program
can be \emph{uniformly delayed} by a $\log \log$ factor on the
time-bound
to run \emph{in $\Lineartime$}.

The theorems of {\bf T} form a \emph{computably enumerable} set, so we
can/do
fix a predicate logic complete automatic theorem prover (such as resolution)
for {\bf T}.  This theorem prover can be time-bounded universally  
simulated --- but, as in the just prior
paragraph, that simulation can be \emph{delayed} by a $\log \log$ factor on 
the time-bound to, then, run \emph{in $\Lineartime$}.
Let
\begin{equation}
{\bf T}\proves_x \fqE
\end{equation}
{\em mean\/} that
\emph{a delayed by such a $\log \log$ factor, 
linear-time computable, time-bounded
universal simulation} of the fixed automatic theorem prover proves $\fqE$
from {\bf T} \emph{within} $x$ steps --- that's \emph{linear-time in}
$(\cardb{\fqE}+\cardb{x})$.

Let $D_x$ be the finite set $(\subseteq \N)$ with canonical index $x$
(see, e.g., \cite{Rog:b:67}).  $x$ codes, for example, {\em both\/} how
to list $D_x$ {\em and\/} how to know
when the listing is done.
Herein, we can and do restrict our canonical indexing of finite sets to
those of sets cardinality $\leq 2$.  We do that in linear-time thus.
Let $0$ be the code of $\emptyset$, and, for set $\set{u,u+v}$,
let the code be $\pair{u,v}+1$. This
coding (suggested by a referee to replace our original one) is linear-time 
codable/decodable (and is 1-1, and, unlike our original, is onto).\footnote{
As an aside: \cite{Cas-Koe:c:09:diff-fair} canonically
codes \emph{any size} finite sets in cubic time \& decodes them
in linear-time.}

\section{Results}\label{results}

   \subsection{Non-Obfuscated Unprovable Programs}\label{main-result}

\begin{thm}\label{thm:exists-g-w-domain-two}
There \emph{exists} $g \in \Lineartime$ and \emph{small}
positive $c,d\in\natnum$ such
that, for \emph{any} $p$, $\cardb{D_{g(p)}} = 2$ and
\emph{there is} a $q \in D_{g(p)}$ for which:
\begin{equation}
\phitm_q = \phitm_p;
\end{equation}
for all $x\in\natnum$,
\begin{equation}\label{one}
\Phitm_q(x) \leq c\cdot(\cardb{p} + \cardb{x} + \Phitm_p(x));
\end{equation}
\begin{equation}\label{two}
\cardb{q} \leq d \cdot \cardb{p};
\end{equation}
yet
\begin{equation}\label{unprovable}
{\bf T} \not\proves \formula{\phitm_q = \phitm_p}.
\end{equation}
\end{thm}

Our proof below of Theorem~\ref{thm:exists-g-w-domain-two},
as will be seen,
makes it \emph{easily transparent} that $\phitm_q = \phitm_p$. Hence,
$q$ is
\emph{not obfuscated}, yet its correctness (at computing $\phitm_p$),
as will
also be seen, is unprovable in {\bf T}. From the \emph{time and program
size
complexity} content of the theorem,
$q$ is nicely \emph{only slightly, linearly} more complex than $p$.
Furthermore, our proof is what is called in
\cite[Page~131]{Roy-Cas:b:94}
\emph{a rubber wall argument}: we set up a rubber wall, i.e., a
{\em potential\/} contradiction
off of which to bounce, so that, were the resultant construction to
veer into
satisfaction of an undesired condition (undesired here is the failure
of~(\ref{unprovable}) above), it bounces off the rubber wall (i.e.,
contradiction) toward our goal, here~(\ref{unprovable}),
instead.\footnote{More discussion on identifying contradictions with
walls, a.k.a.~boundaries, can be found on
\cite[Page~131]{Roy-Cas:b:94}.}

\Prf{of Theorem~\ref{thm:exists-g-w-domain-two}}

\begin{sloppypar}
By \emph{two} applications of \emph{linear-time}:
 krt, if-then-else (these from Section~\ref{complexity-bounded} above),
and $\lam{{\sf E},x}({\bf T}\proves_x{\sf E})$ (this from 
Section~\ref{applied-logics} above),
from \emph{any} $\phitm$-program $p$, one \emph{can algorithmically
find}
in linear-time (in $\cardb{p}$),
programs $e_{1,p}$  and $e_{2,p}$ behaving as follows.
\end{sloppypar}

For each $x$,
\begin{equation}\label{eqn:e1}
	\phitm_{e_{1,p}}(x) = 
\begin{cases}
\phitm_p(x) + 1, &\mbox{if }\T \proves_x
 \ll\phitm_{e_{1,p}} = \phitm_p \gg;\cr
\phitm_p(x), & \mbox{otherwise};\cr
\end{cases}
\end{equation}
and
\begin{equation}\label{eqn:e2}
\phitm_{e_{2,p}}(x) = 
\begin{cases}
0, & \mbox{if }\T \proves_x \ \ll\phitm_{e_{2,p}} = \phitm_p \gg;\cr
\phitm_p(x), & \mbox{otherwise}.\cr
\end{cases}
\end{equation}

Let $g \in \Lineartime$ be such that, for each $p$,  $D_{g(p)} =
\set{e_{1,p}, e_{2,p}}$.  We consider cases regarding $p$ for the choice of
the associated $q \in D_{g(p)}$.

Case~(1). domain$(\phitm_p)$ is infinite.
Suppose for contradiction, for some $x$,
$\T \proves_x \ \ll \phitm_{e_{1,p}} = \phitm_p\gg$.
Since, by assumption, $\T$ {\em does not prove false such sentences},
$\phitm_{e_{1,p}} = \phitm_p$, and by~(\ref{eqn:e1}) above, for all
$x' \geq x$,
$\phitm_{e_{1,p}}(x')$ {\em also\/} $ = \phitm_p(x') + 1$, but, since
domain$(\phitm_p)$ is infinite, we have a contradiction.
Choose $q = e_{1,p}$. Then, trivially, again by~(\ref{eqn:e1}),
$\phitm_{q} = \phitm_p$, but $\T$ does not prove it.

Case~(2). domain$(\phitm_p)$ is finite.
Suppose for contradiction, for some $x$,
 $\T \proves_x \ \ll \phitm_{e_{2,p}} = \phitm_{p} \gg$.
Since, by assumption, $\T$ {\em does not prove false such sentences},
$\phitm_{e_{2,p}} = \phitm_p$, and by~(\ref{eqn:e2}) above, for all $x'
\geq
x$,
$\phitm_{e_{2,p}}(x')$ also $ = 0$, making 
domain$(\phitm_{e_{2,p}})$ infinite,
and, hence, domain$(\phitm_p)$ is infinite, a contradiction.
Choose $q = e_{2,p}$. Then, trivially, again by~(\ref{eqn:e2}),
$\phitm_{q} = \phitm_p$, but $\T$ does not prove it.

In each case, by
$\ifthenelseforus$ and $\krt$
being linear-time (hence, at most linear growth) functions,
$\lam{{\sf E},x}({\bf T}\proves_x{\sf E}) \in \Lineartime$, and
by the complexity upper bounds~(\ref{condefficienttoo})
and~(\ref{complexity-task})
(in Section~\ref{complexity-bounded} above)
as well as the assertion
(in Section~\ref{applied-logics} above)
of the linear-time (and, hence,
linear size) cost of substituting numerals into formulas of $\PA$,
we have small positive $c,d$ such that the theorem's time complexity
bound~(\ref{one}) and it's program size bound~(\ref{two}) above each
hold.
\Qed{Theorem~\ref{thm:exists-g-w-domain-two}}

Next is the promised, motivated, concrete example.

Putnam \cite{Put:i:80} notes that the typical inductive definitions of
grammaticality (i.e., well-formedness)
for propositional logic formulas parallel the typical
definitions of truth (under any truth-value assignment to the
propositional
variables) for such formulas, and that the first kind of inductive
definition provides a {\em short and feasible\/} decision program for
grammaticality.\footnote{In computer science these inductive
definitions
would be called recursive and, as program code, can easily be run
iteratively --- for
efficiency.}~~He goes on to say, though, that the other ways of
providing {\em short and feasible\/} inductive definitions of
such grammaticality which also parallel an inductive definition of
truth are so similar as to constitute {\em intrinsic\/}
grammars (and semantics).  Let $p$ be one of these {\em typical\/}
short and fast
decision procedures for propositional calculus grammaticality
expressed naturally and directly as a $\phitm$-program.  Then by
Theorem~\ref{thm:exists-g-w-domain-two}
above and its proof also above, there is an obviously semantically
equivalent $\phitm$-program $q$ only slightly linearly more complex
than
$p$ in size and run time (so it too is short and feasible); $q$ also
provides the same inductive definition of grammaticality as $p$ which,
then,
parallels the truth definition like $p$ does (after all the else part
of $q$
{\em is\/} $p$ and the if-part of $q$ never comes true); but the
unprovability
(in pre-assigned $\T$) of the semantic equivalence of $q$ with $p$
makes
$q$ a bit peculiar as an {\em intrinsic\/} grammar for propositional
logic,
providing a basis to doubt Putnam's assertion.  However, we do note
that
{\em intensionally\/} \cite{Rog:b:67}
$q$ is a bit unlike $p$ --- since it performs
an always false (quick) test $p$ doesn't.

   \subsection{A Constructivity Concern}\label{constructivity}

    It's interesting to ask: can the condition $\cardb{D_{g(p)}} = 2$
in Theorem~\ref{thm:exists-g-w-domain-two} be improved to
$\cardb{D_{g(p)}} = 1$?  If so, it makes sense to replace a singleton
set, $\set{q}$, by just $q$ and use $g(p)=q$ (not the code of
$\set{q}$).  Anyhow, the answer to the question is, No
(see Theorem~\ref{thm:no-g-for-phi} below).
Before we present and prove this theorem, it is useful to have for its
proof the unsurprising lemma (Lemma~\ref{lemma:pa-proves-re}) just
below.\footnote{We bother to prove it since we do not know a citation
for its proof.}

\begin{lem}\label{lemma:pa-proves-re}
If $\phitm_p(x)\converges = y$, then
\begin{equation}
\PA \proves \ \ll \phitm_p(x)\converges = y \gg.
\end{equation}
\end{lem}

\Prf{of Lemma~\ref{lemma:pa-proves-re}}
The relation, in $p,x,y,t$, that holds iff
$\phitm_p(x)\converges = y \mbox{ within } t \mbox{ steps}$, where the
steps are measured by the natural
$\Phitm$, is trivially computable (a.k.a.~recursive)
\cite{Blu:j:67:complexity}.

Suppose $\phitm_p(x)\converges = y$.  Then
there is some $t$ such that $\phitm_p(x)\converges = y \mbox{ within }
t
\mbox{ steps}$. By G\"odel's Lemma \cite{God:j:31,Men:b:09} that
recursive
relations are numeralwise provably-representable in, e.g., $\PA$,
$\PA \proves \formula{\phitm_p(x)\converges = y \mbox{ within } t
\mbox{ steps}}$.  By existential generalization \emph{inside} $\PA$, we
have
$\PA \proves \formula{(\exists t)[\phitm_p(x)\converges = y \mbox{
within } t
\mbox{ steps}]}$. Hence, $\PA \proves \ \ll \phitm_p(x)\converges =
y\gg$.
\Qed{Lemma~\ref{lemma:pa-proves-re}}

The next theorem implies that, in
Theorem~\ref{thm:exists-g-w-domain-two} above, the
condition $\cardb{D_{g(p)}} = 2$
can{\em not\/} be improved to $\cardb{D_{g(p)}} = 1$ (or equivalent as
discussed above). The proof of this next theorem
(Theorem~\ref{thm:no-g-for-phi}) provides positive cases regarding proving
true
program properties in $\PA$.\medskip

\begin{thm}~\label{thm:no-g-for-phi}
It is \emph{not} the case that there
exists computable $g$ such that, for any $p$,
for $q = g(p)$,
\begin{equation}\label{eqn:notproves}
{\bf T} \not\proves \ll\phitm_q=\phitm_p\gg.
\end{equation}
\end{thm}

    \Prf{of Theorem~\ref{thm:no-g-for-phi}}
Suppose for contradiction otherwise.

Suppose $d$ is a $\phitm$-program for $g$, i.e., suppose $\phitm_d =
g$.

Of course $\lam{p,x}[\phitm_{\phitm_d(p)}(x)]$ is partially computable,
and, importantly, this is provable in $\PA$.
We sketch how we know the provability in $\PA$.

For example, one step in showing the provability
is to explicitly construct a $\phitm$ universal
program $u$ so that its detailed correctness is (trivially, albeit
tediously)
provable in $\PA$. In particular,
\begin{equation}
\PA \proves \formula{ (\forall p,x)[\phitm_u(p,x) = \phitm_p(x)] }.
\end{equation}
For $\phitm$, the construction of a relatively efficient, but
{\em time-bounded variant\/} of such a $u$ is outlined in the proof of
\cite[Theorem~3.6]{Roy-Cas:b:94}. This construction can be altered
to remove the time-boundedness and just get a suitable $u$.

Another step would be to spell out a $\phitm$-program $c$ for a
computable function comp$_2$
for computing a $\phitm$-program for the composition of the partial
functions computed by its $\phitm$-program arguments as in
the $m=2$ case of \cite[Lemma~3.10]{Roy-Cas:b:94}
\footnote{\label{fn-compm}
In that Lemma~3.10, we have, in effect,
for the arbitrary $m$ case,
for all $p_0,\ldots, p_m$,
\begin{equation}\label{compm}
\phi_{comp_m(p_0,\ldots,p_m)}(x) =
\phi_{p_0}(\phi_{p_1}(x), \ldots, \phi_{p_m}(x)):
\end{equation}
}
and its proof --- where, again, $\PA$ proves correctness (including
comp$_2 = \phitm_c$ is total).

Relevance of $u$ and $c$:  clearly we have,
\begin{equation}\label{u}
\phitm_{\phitm_d(p)}(x) = \phitm_u(\phitm_d(p),x),
\end{equation}
and the right-hand side of~(\ref{u}) just above is a relevant
composition and,
then, can be further expanded employing $c$.\footnote{
Further below in Section~\ref{subtleties},
we'll consider, among other things, some programming systems with
provability
subtleties regarding universality and/or composition. \emph{This}
composition, though,
we be as in the $m=1$ case of \cite[Lemma~3.10]{Roy-Cas:b:94} 
(see Footnote~\ref{fn-compm} just above).}
With $u$ and $c$, then, we can explicitly compute a $\phitm$-program
for $\lam{p,x}[\phitm_{\phitm_d(p)}(x)]$
and prove it correct in $\PA$.

So, then, by the Constructive
Kleene's Second Recursion Theorem {\em but without the
parameter $p$ as above in Section~\ref{complexity-bounded} above},
we have a
(self-referential) $p_0$ such that, for any $x$,
\begin{equation}\label{eqn-fp}
\phitm_{p_0}(x) = \phitm_{\phitm_d(p_0)}(x).
\end{equation}
Below we'll refer to this parameter-free version of the above
Constructive
Kleene Theorem as $\KRT$.
Then we have a $\phitm$-program $k$ for the above function
$\krt$ {\em again with parameter $p$ completely omitted}, and, with
this $k$ representing in the language of $\PA$ this modified version of
the function $\krt$, $\KRT$ is completely provable in $\PA$.

Hence, by our remarks above about computing and proving
correct a program for $\lam{p,x}[\phitm_{\phitm_d(p)}(x)]$, we can
explicitly compute a $p_0$ as in~(\ref{eqn-fp}) and prove it correct in
$\PA$; we have in particular,
\begin{equation}\label{pa-eqn-fp}
\PA \proves \formula{\phitm_{p_0} = \phitm_{\phitm_d(p_0)}}.
\end{equation}

However, we don't know enough about $g$ (and $d$)
to know whether we can prove $g$'s totality in
$\PA$ --- including by
representing $g$ as $\phitm_d$; fortunately, we won't need that.

We do know (at least outside $\PA$) that $g$ is total (since it's a
consequence
of $g$'s assumed computability).  Hence, we know
(at least outside $\PA$) that $g(p_0)\converges$.  Since, from above,
$d$ is a $\phitm$-program for $g$, we have that
$\phitm_d(p_0)\converges = $
to some explicit numerical value $q_0$.
Therefore, from Lemma~\ref{lemma:pa-proves-re} above,
\begin{equation}\label{eqn:converges}
\PA \proves \ \ll \phitm_d(p_0)\converges = q_0 \gg.
\end{equation}
Hence, by substitution of equals for equals and reflexivity of equals
{\em inside\/} $\PA$,~(\ref{pa-eqn-fp}), and~(\ref{eqn:converges}),
\begin{equation}
\PA \proves \formula{ \phitm_{q_0} = \phitm_{p_0} },
\end{equation}
a contradiction to our beginning assumption --- since $\T$ extends
$\PA$.
    \Qed{Theorem~\ref{thm:no-g-for-phi}}

   \subsection{Subtleties}\label{subtleties}

So far we have considered the natural, deterministic
complexity theory relevant,
acceptable system, $\phitm$.  After we obtained
Theorem~\ref{thm:no-g-for-phi}
just above
--- which shows a condition in Theorem~\ref{thm:exists-g-w-domain-two}
further above
(in Section~\ref{main-result}) couldn't be improved,
we wondered if there were some (possibly not quite so natural but,
perhaps,
still acceptable) systems $\psi$ for which we don't have the just above
Theorem~\ref{thm:no-g-for-phi}.  We initially obtained the first part
of the
next theorem (Theorem~\ref{thm:exists-a-g})
which provides such a $\psi$, {\em but\/}
we didn't, then,
know whether our $\psi$ was acceptable.  We subsequently obtained
Theorem~\ref{thm:exists-a-g}'s furthermore clause providing our
$\psi$'s
acceptability
{\em together with a surprise we didn't expect}.  We explain the
surprise
after the statement of
Theorem~\ref{thm:exists-a-g} and before its proof.

\begin{thm}\label{thm:exists-a-g}
There is a programming system $\psi$
and a computable $g$ such that, for all $p$, $\psi_{g(p)} = \psi_p$,
yet, for $q = g(p)$,
$ \T \not \proves \ \ll\psi_q=\psi_p\gg$.

Furthermore, $\psi$ is acceptable, and, surprisingly,
\begin{equation}\label{surprise}
(\forall p)[\psi_p = \phitm_p].
\end{equation}
\end{thm}

\noindent How can~(\ref{surprise}) be true --- in the light of the rest of the
just
above theorem (Theorem~\ref{thm:exists-a-g})?  It seems to contradict
Theorem~\ref{thm:no-g-for-phi} further above.
The answer is that, in the proof just below of the just above theorem
(Theorem~\ref{thm:exists-a-g}), the needed $\psi$ is, in effect,
{\em defined by\/} an unusual
$\phitm$-program $e$ in~(\ref{eqn:obfuc-e-psi}, \ref{eqn:obfuc-e})
below,
and, {\em in the language of\/} $\PA$, for formulating
(un)provability {\em about\/} $\psi$ in $\T$,
$\psi$ is, of course, {\em represented by its defining}~$e$.\footnote{
An original source for unusual representations in arithmetic
(as is our $e$) is \cite{Fef:j:60}.}  $\phitm$ itself, on the other
hand,
can be and is understood to be naturally (not unusually)
represented in the language of
$\PA$.\footnote{See the informal discussion about the notation
$\formula{\sf E}$ in Section~\ref{applied-logics} above,
where, in effect, the particular
example $\formula{\phitm_p \mbox{ is total}}$ is employed.}

To aid us in some proofs below, including that of the above
Theorem~\ref{thm:exists-a-g}, we present the following lemma
(Lemma~\ref{lemma:constructive-padding}), where the recursion theorem
part of its proof is from H.~Friedman~\cite{Fri:m:74}.
\begin{lem}[{\bf T}-Provable Padding-Once]
					\label{lemma:constructive-padding}
Suppose $\alpha$ is any acceptable programming system such that
$\T$ proves $\alpha$'s acceptability.

Then, there is a total computable function
$g$ such that for any $p$, $g(p) \neq p$, but $\alpha_{g(p)} =
\alpha_p$.\footnote{
Of course, a more general, constructive {\em infinite\/} padding holds
\cite{Mac-You:b:78}, and we need a version of that further below.}

Furthermore, this padding-once result is, then, expressible and
provable in $\T$.

\end{lem}

\Prf{of Lemma~\ref{lemma:constructive-padding}}
Assume the hypothesis, i.e., that $\T$ proves $\alpha$'s acceptability.

Then $\T$ proves Kleene's S-m-n Theorem, so we obtain
that $\T$ proves
the Parameterized Second Kleene Recursion Theorem
(as above in Section~\ref{complexity-bounded},
but with witnessing functions
not necessarily in $\Lineartime$).

Then, from this Kleene Theorem,
we have a computable function $f$
such that, for each $p,x$,
\begin{equation}\label{eqn:padding-f}
\alpha_{f(p)}(x) =\begin{cases}
	\alpha_p(x), & \mbox{if }f(p) \neq p;\cr
	\alpha_{p+1}(x), & \mbox{if }f(p) = p.\cr
\end{cases}
\end{equation}
Then, let $g$ be defined as follows.
\begin{equation}\label{eqn:padding-g}
g(p) = \begin{cases}
	f(p), & \mbox{if }f(p) \neq p; \cr
	p+1, & \mbox{if }f(p) = p.\cr
\end{cases}
\end{equation}
We consider two cases.
\begin{description}
\item[Case one] $f(p) \neq p$.
Then, from~(\ref{eqn:padding-f}),
$\alpha_{f(p)} = \alpha_p$, and, from~(\ref{eqn:padding-g}),
$g(p) = f(p) \neq p$.

\item[Case two] $f(p) = p$.
Then, from~(\ref{eqn:padding-f}),
$\alpha_{f(p)} = \alpha_{p+1}$, which, by Case~two, $ = \alpha_p$.
From~(\ref{eqn:padding-g}),
$g(p) = p+1 \neq p$.
\end{description}
By this case-analysis, $g$ satisfies Padding-Once.  The above is
so simple as to be provable in $\T$ --- as needed.
\Qed{Lemma~\ref{lemma:constructive-padding}}

\Prf{of Theorem~\ref{thm:exists-a-g}}
By Kleene's second recursion theorem (again without parameter), there
is
a (self-referential) $\phitm$-program $e$ and an associated $\psi$
both such that, for each $p,x$,
\begin{equation}\label{eqn:obfuc-e-psi}
\psi_p(x) \eqdef \phitm_e(\pair{p,x}), \mbox{ which } =
\end{equation}
\begin{equation}\label{eqn:obfuc-e}
\begin{cases}
	p, & \mbox{if }\T \proves_x
	    \ll(\exists q,r \mid q \neq r)[\psi_q = \psi_r]\gg;\cr
	\phitm_{p}(x), & \mbox{otherwise}.
\end{cases}
\end{equation}
N.B. The mentions of $\psi$ in~(\ref{eqn:obfuc-e}) just above with
variable
subscripts $q,r$
should be understood, employing $\psi$'s
definition~(\ref{eqn:obfuc-e-psi})
above, to be
$\phitm_e(\pair{q,\cdot}), \phitm_e(\pair{r,\cdot})$, respectively.

\begin{clm}\label{claim:t-not-prove-eq} $\T \not\proves \ \ll(\exists
q, r \mid q \neq r)[\psi_q = \psi_r]\gg$.
\end{clm}

\Prf{of Claim~\ref{claim:t-not-prove-eq}}
Suppose for contradiction otherwise.

Then there
exists an $x_0$ such that $\T \proves_{x_0} \ \ll(\exists q, r \mid q
\neq r)[\psi_q = \psi_r]\gg$.  However, then,
by~(\ref{eqn:obfuc-e-psi},
\ref{eqn:obfuc-e}) above, we have $(\forall
p)[\psi_p(x_0)\converges = p]$, and thus $(\forall p, q \mid p \neq
q)[\psi_p(x_0) \neq \psi_q(x_0)]$; therefore, $\T$ has proven a
sentence of first order arithmetic which is false in the standard
model, a contradiction.
\Qed{Claim~\ref{claim:t-not-prove-eq}}

\begin{clm}\label{claim:psi-eq-phi}
$(\forall p)[\psi_p = \phitm_p]$; hence, $\psi$ is acceptable.
\end{clm}

\Prf{of Claim~\ref{claim:psi-eq-phi}}
By Claim~\ref{claim:t-not-prove-eq},
the first clause in (\ref{eqn:obfuc-e}) above is false for each $p,x$.
Therefore, by~(\ref{eqn:obfuc-e-psi}, \ref{eqn:obfuc-e}) above,
$(\forall p,x)[\phitm_e(\pair{p,x}) = \phitm_p(x)]$; hence,
$(\forall p)[\psi_p = \phitm_p]$ --- making $\psi$ acceptable too.
\Qed{Claim~\ref{claim:psi-eq-phi}}

\begin{clm}\label{claim:comp-g-exists} There is a computable $g$ such
that,
for all $p$, $g(p) \neq p$,  $\psi_{g(p)} = \psi_p$, and, for $q =
g(p)$, $\T \not \proves \ \ll\psi_q = \psi_p\gg.$
\end{clm}

\Prf{of Claim~\ref{claim:comp-g-exists}}
The acceptability of $\phitm$
is  provable in $\PA$, hence, in $\T$.
By Lemma~\ref{lemma:constructive-padding}, there exists a
computable $g$ such that $(\forall p)[g(p) \neq p \AND \phitm_{g(p)} =
\phitm_p]$.  By Claim~\ref{claim:psi-eq-phi}, $\psi = \phitm$; thus,
for
this \emph{same} $g$, $(\forall p) [\psi_{g(p)} = \psi_p]$.

Suppose arbitrary $p$ is given.  Let $q = g(p)$.
Suppose for contradiction $\T \proves \ \ll\psi_q = \psi_p\gg$.
Clearly by G\"odel's Lemma (employed in the proof of
Lemma~\ref{lemma:pa-proves-re} above), $\PA \proves \formula{q \neq
p}$.
Then, by this and existential generalization in $\T$,
$\T \proves \ \ll(\exists q, r \mid q \neq r)[\psi_q = \psi_r]\gg$,
a contradiction to Claim~\ref{claim:t-not-prove-eq} above.
\Qed{Claim~\ref{claim:comp-g-exists}}

\Qed{Theorem~\ref{thm:exists-a-g}}


For the next three corollaries
(Corollaries~\ref{cor:t-not-prove-psi-eq-phi},~\ref{cor:not-krt-and-u},
and~\ref{cor:t-not-prove-eps}), the mentioned $\psi$ is that from
Theorem~\ref{thm:exists-a-g} and its proof,
including~(\ref{eqn:obfuc-e-psi}, \ref{eqn:obfuc-e}) above.

\begin{cor}\label{cor:t-not-prove-psi-eq-phi} $\psi = \phitm$, but $\T
\not\proves \ \ll \psi = \phitm \gg$.\end{cor}

\Prf{of Corollary~\ref{cor:t-not-prove-psi-eq-phi}}
$\psi = \phitm$ is from Theorem~\ref{thm:exists-a-g} above.
Suppose for contradiction
$\T \proves \ \ll \psi = \phitm \gg$.

Then, from this and the proof of Theorem~\ref{thm:no-g-for-phi} above,
one obtains a Theorem~\ref{thm:no-g-for-phi} {\em but\/}
with $\psi$ replacing $\phitm$.
This contradicts Theorem~\ref{thm:exists-a-g} above (which is also
about $\psi$).
\Qed{Corollary~\ref{cor:t-not-prove-psi-eq-phi}}

To understand the corollary (Corollary~\ref{cor:not-krt-and-u}) and its
proof just below,
it may be useful to review the roles of $\phitm$-programs $u,c,k$
in the proof of Theorem~\ref{thm:no-g-for-phi} above. This corollary
says there can be {\em no\/} analog of {\em all three\/} of these
programs
for $\psi$ (in place of $\phitm$).

\begin{cor}\label{cor:not-krt-and-u} There are \emph{no} $u,c,k$ such
that {\em simultaneously:}
\begin{equation}
\T \proves \ \ll u \mbox{ is a witness to universality in }  \psi\gg,
\end{equation}
\begin{equation}
\T \proves \ \ll c \mbox{ is a witness to composition in } \psi\gg\footnote{
\emph{This} composition is the $m=2$ case of
\cite[Lemma~3.10]{Roy-Cas:b:94} (see Footnote~\ref{fn-compm} further above).},
\mbox{\&}
\end{equation}
\begin{equation}\label{third-one}
\T \proves \ \ll k \mbox{ is a witness to } \KRT \mbox{ in } \psi\gg.
\end{equation}
\end{cor}

\Prf{of Corollary~\ref{cor:not-krt-and-u}}
Suppose for contradiction otherwise.  Then, enough is provable in $\T$
{\em about\/} $\psi$
to make Theorem~\ref{thm:no-g-for-phi} above
{\em also\/} provable for $\psi$ --- in place of $\phitm$.
This contradicts Theorem~\ref{thm:exists-a-g} about $\psi$.
\Qed{Corollary~\ref{cor:not-krt-and-u}}

Regarding Corollary~\ref{cor:not-krt-and-u} just above, it is well known that,
from \cite{Mac-You:b:78,Mar:j:93}, Kleene's S-m-n can be constructed
out
of a program $c$ for composition
and, then, Kleene's proof of $\KRT$
can be done from S-m-n; so, it might appear that~(\ref{third-one})
just above could be eliminated. This is actually open.  The reason is
that,
while each of these just mentioned
constructions requires some easily existing auxiliary $\psi$-programs,
for Corollary~\ref{cor:not-krt-and-u} we'd ostensibly also need these
auxiliary $\psi$-programs to be $\T$-provably correct.\footnote{Machtey
and Young's construction
\cite{Mac-You:b:78,Mar:j:93} of an S-m-n function out of a
composition function, for example, employs auxiliary $\psi$-programs
$q_0,q_1$ such that
\begin{equation}
\psi_{q_0} = \lam{z}\pair{0,z}; \mbox{ and }
\psi_{q_1} = \lam{\pair{y,z}}\pair{y+1,z}.
\end{equation}
Marcoux's more efficient solution \cite{Mar:j:93} employs
three such auxiliary $\psi$-programs.}
The $\T$-provable correctness is the hard part.

\begin{cor}\label{cor:t-not-prove-eps} $\T \not \proves \ \ll\psi $ is
acceptable$ \gg$.
\end{cor}

\Prf{of Corollary~\ref{cor:t-not-prove-eps}}
Assume for contradiction otherwise.  Then,
by Lemma~\ref{lemma:constructive-padding} above,
$\T \proves \ \ll(\exists p)[\psi_{p} \mbox{ is
total} \AND
(\forall q)(\exists r = \psi_{p}(q) \mid r \neq q)[\psi_q =
\psi_r]\gg$.

From this we have,
$\T \proves \ \ll(\exists q, r \mid q \neq r)[\psi_q =
\psi_r]\gg$, a contradiction to Claim~\ref{claim:t-not-prove-eq}
above.
\Qed{Corollary~\ref{cor:t-not-prove-eps}}

\subsubsection{Subtleties About Proving Universality}\label{univ}

The next two theorems herein (Theorems~\ref{thm:t-prove-one-u}
and~\ref{thm:t-not-prove-u-in-theta}) provide two more acceptable
programming systems, $\eta,\theta$,
respectively, each defined (as was $\psi$ above)
by respective, unusual $\phitm$-programs.
The first of these theorems (Theorem~\ref{thm:t-prove-one-u})
provides a surprise, part positive, part negative, regarding
proving {\em in} $\T$ that \emph{universality} holds for $\eta$. The
contrast between these last two theorems is also interesting.   Of
course,
since each of $\eta,\theta$ is acceptable, universality holds for each
of them
(at least outside $\T$).

 Below, for partial functions $\xi$, $\rho(\xi)$ denotes the
 {\em range\/} of $\xi$.

\begin{thm}\label{thm:t-prove-one-u}
There exists an acceptable programming system $\eta$ and an $e$
such that $\PA \proves \ \ll e$ is universal for $\eta\gg$, yet,
surprisingly, for each $p$,
\begin{equation}
\mbox{If } \T \proves \ \ll p \mbox{ is universal for } \eta\gg,
\mbox{ then } p = e.
\end{equation}
Of course, in $\eta$, there are infinitely many universal programs, but
\emph{exactly one} provably so in $\T$. 
Furthermore, $\eta$ turns out to be $\phi$.
\end{thm}

\Prf{of Theorem~\ref{thm:t-prove-one-u}}
The Kleene Second Recursion Theorem
provides a $\phitm$-program $e$ and an associated
$\eta$ both such that, for each $p,x$,
\begin{equation}\label{eqn:obfuc-e-eta}
\eta_p(x) \eqdef \phitm_{e}(\pair{p,x}), \mbox{ which } =
\end{equation}
\begin{equation}\label{eqn:obfuc-e-prime}
\begin{cases} p,   & \mbox{if }[p \neq e \AND \T \proves_x
	      \ll(\exists q,r \mid r \neq  q) [\eta_q =
	     \eta_r]\gg];\cr
\phitm_p(x), & \mbox{otherwise}.
\end{cases}
\end{equation}\smallskip

\noindent Of course, since $\KRT$ for $\phitm$ is constructively provable in
$\PA$,
we can get the numeral for $e$ inside $\PA$ as well as the universally
quantified
equation just above for the value of $\phitm_{e}(\pair{p,x})$ by
cases.

\begin{clm}\label{claim:t-not-prove-for-eta}$\T \not \proves \ \ll
(\exists q, r \mid q \neq r)[\eta_q = \eta_r] \gg$.
\end{clm}

\Prf{of Claim~\ref{claim:t-not-prove-for-eta}}
Assume for contradiction that 
\begin{equation}\label{not-false}
\T \proves \formula{(\exists q, r \mid q \neq r)[\eta_q = \eta_r]}.  
\end{equation}
Let $x_0$ be the minimum number of steps in any
such proof. Since $\T$ does not prove false sentences of $\PA$ and 
(\ref{not-false}), we have
\begin{equation}\label{eqn:true}
(\exists q, r \mid q \neq r)[\eta_q = \eta_r].
\end{equation}
Then, for $f(p) = \phitm_{e}(\pair{p,x_0})$,
$\rho(f) \supseteq (\natnum-\set{e})$.
Clearly, $\eta_{e} = \phitm_{e}$, and $\eta_{e}$ has infinite
range. Furthermore, $(\forall p \neq e) (\forall x \geq x_0)[\eta_p(x)
=
p]$; therefore, $(\forall p \neq e)[\eta_p$ has finite range$]$, and,
thus, there is no $\eta$-program whose code number is not $e$ whose
computed
partial function is equal to $\eta_{e}$.
Furthermore, $(\forall p,q \mid p \neq e \AND p \neq q \AND q
\neq e) [\eta_p(x_0) = p \AND \eta_q(x_0) = q]$, thus there are no
two distinct programs that compute the same partial function,
a contradiction
to~(\ref{eqn:true}).
\Qed{Claim~\ref{claim:t-not-prove-for-eta}}

\begin{clm}\label{claim:eta-eq-phi-2}
$(\forall p, x)[\eta_p(x) = \phitm_p(x)]$;
hence, $\eta$ is acceptable.  \end{clm}

\Prf{of Claim~\ref{claim:eta-eq-phi-2}}
By Claim~\ref{claim:t-not-prove-for-eta}, the
if clause of~(\ref{eqn:obfuc-e-prime}) is always false,
hence, by~(\ref{eqn:obfuc-e-eta}, \ref{eqn:obfuc-e-prime}),
$(\forall p, x)[\eta_p(x) = \phitm_{e}(\pair{p,x}) = \phitm_p(x)]$.
\Qed{Claim~\ref{claim:eta-eq-phi-2}}

\begin{clm}\label{claim:t-not-prove-equiv-in-eta-2} There does not
exist $p, q$ such that $p \not= q$, and $\T \proves \ \ll \eta_p =
\eta_q \gg$. \end{clm}

\Prf{of Claim~\ref{claim:t-not-prove-equiv-in-eta-2}}
Suppose for contradiction otherwise.  Then, by G\"odel's Lemma followed
by
existential generalization,
the latter in $\T$, we obtain a contradiction to
Claim~\ref{claim:t-not-prove-for-eta}.
\Qed{Claim~\ref{claim:t-not-prove-equiv-in-eta-2}}

\begin{clm}\label{claim:t-prove-e-u-in-eta} $\PA \proves \ \ll e$ is
universal for $\eta\gg$.
\end{clm}

\Prf{of Claim~\ref{claim:t-prove-e-u-in-eta}}
We need {\em not\/} prove in $\PA$ that $\eta$ is a programming system
for
the 1-argument partial computable functions.
Instead, it suffices for us to argue only that
\begin{equation}\label{eqn:eta.univ}
 \PA \proves \ \formula{(\forall p, x)[\eta_e(\langle p,x \rangle) = \eta_p(x)]}.
\end{equation}
Then, from~(\ref{eqn:obfuc-e-eta}) above, the definition of $\eta$ by
$e$
in the $\phitm$-system, applied to each side of~(\ref{eqn:eta.univ}),
it, then, suffices to show that
\begin{equation}\label{eqn:phitm.univ}
\PA \proves \formula{(\forall p, x)[
\phitm_{e}(\pair{e,\pair{p,x}}) = \phitm_{e}(\pair{p,x})]}.
\end{equation}
By the otherwise clause of~(\ref{eqn:obfuc-e-eta},
\ref{eqn:obfuc-e-prime})
above, applied to $\formula{\phitm_{e}(\pair{e,\pair{p,x}})}$, where
$p,x$ are
variables (not numerals), we get its provable
in $\PA$ value to be $\formula{\phitm_{e}(\pair{p,x})}$ --- again with
$p,x$
variables. This together with
universal generalization {\em inside\/} $\PA$ on the variables $p,x$,
verifies {\em in\/}
$\PA$ the sufficient~(\ref{eqn:phitm.univ}) just above.
\Qed{Claim~\ref{claim:t-prove-e-u-in-eta}}

\begin{clm}\label{claim:t-prove-one-u}
For all $p \not= e$, $\T
\not\proves \ \ll p$ is universal in $\eta\gg$.\end{clm}

\Prf{of Claim~\ref{claim:t-prove-one-u}}
Immediate from Claims
\ref{claim:t-not-prove-equiv-in-eta-2} and
\ref{claim:t-prove-e-u-in-eta}.
\Qed{Claim~\ref{claim:t-prove-one-u}}

\Qed{Theorem~\ref{thm:t-prove-one-u}}

\begin{thm}\label{thm:t-not-prove-u-in-theta}
There exists an acceptable programming system  $\theta$
such that, for each $u$,
\begin{equation}
\T \not \proves \ \ll u \mbox{ is universal in } \theta \gg.
\end{equation}
Of course, in $\theta$, there are infinitely many universal programs,
but
\emph{none} are provably so in $\T$. Furthermore, $\theta$ turns out to be
$\phi$.
\end{thm}

\Prf{of Theorem~\ref{thm:t-not-prove-u-in-theta}}
Kleene's Recursion Theorem
provides a $\phitm$-program $e$ and an associated
$\theta$ both such that, for each $p,x$,
\begin{equation}\label{eqn:e-for-theta}
\theta_p(x) \eqdef \phitm_{e}(\pair{p,x}) \mbox{ which } =
\end{equation}
\begin{equation}\label{eqn:e-for-theta-prime}
\begin{cases} p, & \mbox{if }\T \proves_x
    \ll(\exists u)[u $ is universal in $ \theta]\gg; \cr \phitm_p(x),
   & \mbox{otherwise}.
\end{cases}
\end{equation}
Assume for contradiction that $\T \proves \ \ll(\exists u)[u $ is
universal in $ \theta]\gg$.  Then, since $\T$ does not prove false
things
of this sort, universality holds in $\theta$.

Let $x_0$ be the smallest
number of steps in any proof as is assumed just above to exist.

Then, since $(\forall p, x \mid x \geq x_0)[\theta_p(x)=p]$, we have,
$(\forall p)[|\rho(\theta_p)| \leq 1+x_0]$. By contrast,
 $\rho(\theta) = \N$.
Then there is no $p$ such that $\rho(\theta_p) = \rho(\theta)$;
therefore,
there cannot be any universal programs for $\theta$, a contradiction.

Therefore, $\T \not \proves \ \ll(\exists u)[u $ is universal in $
\theta]\gg$, and, thus, $(\forall p, x)[\theta_p(x) = \phitm_p(x)]$.
This makes $\theta$ acceptable.

Furthermore, it is not the case that $(\exists u)[\T \proves \ \ll u $
is universal in $ \theta\gg]$, since, if $\T$ proved such a thing, it
would
immediately follow from Existential Generalization {\em in\/} $\T$
that $\ll(\exists u)[u $ is universal in $ \theta]\gg$ is provable in
$\T$,
which has already been shown not to be provable by $\T$.
\Qed{Theorem~\ref{thm:t-not-prove-u-in-theta}}

We expect that analogs of Theorems~\ref{thm:t-prove-one-u}
and~\ref{thm:t-not-prove-u-in-theta} just above
can be obtained for other properties besides universality. 
In the next section we have an analog for composition.

\subsubsection{Provable Composition with
Unprovable Program Equivalence}\label{comp}

The main result of this section 
(Theorem~\ref{thm:t-prove-composition-not-equality} below):
there is an acceptable programming system (equivalent
to $\phitm$) such that $\T$ can prove there exists a specific program
which
witnesses composition in that system, but $\T$ is \emph{still} unable
to prove
that there exist two distinct, equivalent programs in that system.
It is the last and hardest to prove result in the present paper.

In the proof of this theorem,
we require provable in $\PA$ \emph{infinite} padding for just
the $\phitm$-system --- a stronger form of padding than provided by
Lemma~\ref{lemma:constructive-padding}, though this \emph{infinite}
padding function is only valid for $\phitm$, and thus would not have
been usable for Corollary~\ref{cor:t-not-prove-eps}.
As such, we introduce the following function $\pad$, which uses the
concepts from
\cite[Chapter~3 but as modified in the associated Errata]{Roy-Cas:b:94}
of normal and abnormal code numbers of $\phitm$-programs.  Herein we
briefly discuss these concepts.
A program in the $\phitm$-system is defined
as a \emph{non}-empty
sequence of instructions for a $k$-tape Turing machine, where all
these instructions have the same $k$.
If a given number is not directly the code for such
a sequence, it is defined to be an \emph{abnormal} code; otherwise, it
is a \emph{normal} code. The process of checking a given number to
see if it is normal or abnormal is computable.\footnote{In fact, it is
\emph{linear-time} checkable although nothing in the proofs herein
makes use of that fact.}
By convention, abnormal codes are treated as each encoding the same
Turing
machine, and thus all compute the same function. It is a consequence of
the
encoding used for the $\phitm$-system that all normal codes are
divisible by eight; thus, there are infinitely
many even abnormal codes --- for instance, anything divisible by two
but not eight is an even abnormal code.

\begin{equation}\label{eqn:obfuc-comp-pad}
\pad(p) \eqdef
\begin{cases}
\mbox{the next even abnormal code}, & \mbox{if }p\mbox{ is an abnormal
code};\cr
p\mbox{ with the last instruction repeated}, & \mbox{if }p\mbox{ is a
normal code}.
\end{cases}
\end{equation}

\begin{clm}\label{claim:pad-works}
For any input $p$, the output of $\pad(p)$ is the code number of an
even
program $q$ such that $q > p$, $\phitm_q = \phitm_p$, and $\PA$ proves
this.
\end{clm}
\Prf{of Claim~\ref{claim:pad-works}}
As above,
determining whether a given number is a normal code is
algorithmic. Furthermore, as there are infinitely many even abnormal
codes, finding the next such is straightforward.  Thus, the first
clause in
the definition of $\pad$ is computably checkable and, if true, $\pad$
outputs an even number greater than the input, such that both are
abnormal codes. As noted above,
each abnormal code is defined (for the $\phitm$-system)
to compute the same function as each other abnormal code,
and, hence, the part of the claim not about $\PA$ holds for such $p$.

If the second clause holds, then converting $p$ into a coded sequence
of
instructions is simple, repeating the last instruction is also simple,
and
converting that sequence back into a code number is once again simple.
By \cite{Roy-Cas:b:94}, when there are multiple instructions that apply
in a
given state, the \emph{first} one is the one that applies, thus
repeating an instruction can have no impact on the (possibly partial)
function computed by a $\phitm$ program; thus, $\phitm_{\pad(p)} =
\phitm_p$ in this case. Furthermore, as noted above,
all normal instruction codes \emph{must} be even;
thus, $\pad(p)$ is even in this case as well.

Therefore, the part of the
claim not about $\PA$ holds. The proof so far is so
simple that it can be carried out in $\PA$ straightforwardly albeit
tediously.
\Qed{Claim~\ref{claim:pad-works}}

In our proof of the next theorem
(Theorem~\ref{thm:t-prove-composition-not-equality}) it is convenient
to
employ the following lemma (Lemma~\ref{mixed-rt}), a new recursion
theorem
which mixes an $n$-ary (non-parameter) version of the original
Kleene Recursion Theorem \cite[Page~214]{Rog:b:67} with the Delayed
Recursion Theorem \cite[Theorem~1]{Cas:j:74}.\footnote{In
\cite{Cas:j:74}
the Delayed Recursion Theorem is used to prove the Operator Recursion
Theorem \cite{Cas:j:74,Cas:j:94:self}.
In \cite{Cas:j:91}
the proof of its Remark~1 employs a \emph{subrecursive}
Delayed Recursion Theorem.}

\begin{lem}[A Mixed Recursion Theorem]\label{mixed-rt}
Suppose $n>0$.  Suppose $\xi_1,\ldots,\xi_n,\xi$ are partial
computable.

Then there are $e_1,\ldots,e_n,c$ such that $\phitm_c$ is total, and,
for each $i$ with $1\leq i \leq n$, for all $x,y$,
\begin{equation}\label{nkrt-part}
  \phitm_{e_i}(y) =
\xi_{i}({e_1, \ldots, e_n, c, y}),
\end{equation}
and
\begin{equation}\label{delayed-part}
\phitm_{\phitm_c(x)}(y) =
\xi(e_1, \ldots, e_n, c, x, y).
\end{equation}
\end{lem}\smallskip

\noindent In Lemma~\ref{mixed-rt} just above, (\ref{nkrt-part}) expresses the
$n$-ary Kleene Recursion Theorem part, and (\ref{delayed-part})
expresses the Delayed Recursion Theorem part.

\begin{thm}\label{thm:t-prove-composition-not-equality}
There is an acceptable programming system $\zeta$ and a $w$ such that $
\PA
\proves \ \ll w$ is a witness to composition in $\zeta\gg$, yet $ \T
\not\proves \ \ll(\exists r, t \not= r)[\zeta_r = \zeta_t]\gg$.
Furthermore, $\zeta$ turns out to be $\phi$.
\end{thm}

\Prf{of Theorem~\ref{thm:t-prove-composition-not-equality}}
We apply the $n=3$ case of
Lemma~\ref{mixed-rt} just above to obtain programs $e, c, w', w$
behaving
as below in
(\ref{eqn:obfuc-comp-psi}, 
\ref{eqn:obfuc-comp-c}, 
\ref{eqn:obfuc-comp-w'},
\ref{eqn:obfuc-comp-w}), respectively.\footnote{
Our application of Lemma~\ref{mixed-rt} here does not and does not
need to make full use of
all the self/other reference available in this lemma.}
The $\zeta$ we need for the theorem is defined in terms of this $e$
thus.  For each $p,x$,
\begin{equation}\label{psi}
\zeta_{p}(x) \eqdef \phitm_e(p, x).
\end{equation}
For convenience below, in describing the behavior of
$e, c, w', w$, in many places
we'll write this $\zeta$ instead of $\phitm_e$.

In the following formula we let $\myprime(n)$
be the $n$th prime, where $\myprime(0) = 2$, $\myprime(1) = 3$,
$\ldots$~. Importantly to the combinatorics of the diagonalization below
in this proof, any odd prime raised to a power is odd, and we noted
before the statement of the present theorem being proved 
(Theorem~\ref{thm:t-prove-composition-not-equality}) that our particular
infinite padding function $\pad$ always outputs even numbers.

For each $p,x$,
\begin{equation}\label{eqn:obfuc-comp-psi}
\phitm_e(p, x) =
\begin{cases}
\phitm_p(x), & \mbox{if }p = w \Or (\exists y < p)[\phitm_w(y) = p]
\Or\cr& \Quads{1}\T \not\proves_x \ \ll (\exists r, t \not= r)[\zeta_r
= \zeta_t]\gg;\cr
\myprime(p+1)^x, & \mbox{if }p \not= w \AND (\forall y <
p)[\phitm_w(y)\converges \not= p] \AND\cr& \Quads{1}\T \proves_x \ \ll
(\exists r, t \not= r)[\zeta_r = \zeta_t]\gg;\cr
\diverges, & \mbox{otherwise}.
\end{cases}
\end{equation}
$\phitm$-program $c$, spelled out just below, outputs a $\phitm$-program
which computes the composition of two input $\zeta$-programs.  It is employed
by $\phitm$-program $w'$ further below.  

For each $p,q$,
\begin{equation}\label{eqn:obfuc-comp-c}
\phitm_{\phitm_c(\pair{p, q})}(x) = \zeta_p(\zeta_q(x)).
\end{equation}
$\phitm$-program $w'$, spelled out next, rewrites each output of 
$\phitm$-program $c$ above so that the output of $w'$
computes the same partial
function as this output of $c$, but has combinatorially useful
numeric properties regarding evenness.

For each $p,q$,
\begin{equation}\label{eqn:obfuc-comp-w'}
\phitm_{w'}(\pair{p, q}) =
\begin{cases}
$the first value $v$, if any, found by iteratively
applying$\cr\Quads{1}\pad$ to $\phitm_c(\pair{p, q})$ such that: $v
\not= w$, $v$ is even,$\cr\Quads{1}
v > \pair{p, q}$ (hence, $v > p,q$ by Lemma~\ref{lem-pairing}),
$\cr\Quads{1}(\forall x < \pair{p, q})[\phitm_{w'}(x)\converges]$,
$\cr\Quads{1}[$if $\pair{p, q} > 0$, then $v > \phitm_{w'}(\pair{p,
q}-1) \mbox{ which } \converges]$, $\cr
\Quads{1}$and, for $\pair{r, s} = v, \cr 
\Quads{1}(\neg\exists x < \pair{p, q})[\phitm_{w'}(x) = r \Or
\phitm_{w'}(x) = s];\cr
\diverges, \mbox{if no such $v$ exists.}\cr
\end{cases}
\end{equation}

\noindent As we'll see, $\phitm$-program $w$, spelled out next, 
is such that $\phitm_w = \zeta_w$~(Claim~\ref{claim:psi_w_is_phitm_w} below),
\emph{and} it computes the $m=1$ case control structure of composition 
from \cite[Lemma~3.10]{Roy-Cas:b:94} (see Footnote~\ref{fn-compm} above) ---
but for the 
$\zeta$-system~(Claim~\ref{claim:PA-proves-c-witness-composition}).
This $w$ 
rewrites the output of $w'$ above so that the resultant $\zeta$-system
composition will have strong associativity properties
\emph{at the $\zeta$-program code number level}: let $\comp_1 = \zeta_w$; then
for all $\zeta$-programs $a,b,c$,
the $\zeta$-program number $\comp_1(\comp_1(a,b),c)$ 
\emph{will be the same $\zeta$-program number as} $\comp_1(a,\comp_1(b,c))$.
This strong property is, in effect, further developed in 
Claims~\ref{claim:psi_chain_x_is_psi_x} 
through~\ref{claim:partial-unchain-is-ok}, 
these claims put limits on the nature of \emph{un}equal
$\zeta$-program numbers \emph{in}
$\rho(\comp_1)$ --- so they cannot interfere with
the unprovability part of Theorem~\ref{thm:t-prove-composition-not-equality},
and, then, they are used in proving the difficult to prove
Claim~\ref{claim:T-proves-equality-implies-infinitely-not-equal}. 
This latter claim provides most of the desired unprovability --- with
Claim~\ref{claim:T-not-proves-equality} finishing it off.

For each $p,q$,
\begin{equation}\label{eqn:obfuc-comp-w}
\phitm_w(\pair{p, q}) =
\begin{cases}
\phitm_w(\pair{\phitm_w(\pair{p, r}), s}), & \mbox{if }(\forall
\pair{r,s} < q)[\phitm_w(\pair{r,s})\converges] \AND\cr & \Quads{1}
(\exists \pair{r,s} < q)[\phitm_w(\pair{r,s}) = q],\cr & \Quads{1}
\mbox{then select minimum such $s$}\cr
& \Quads{1}\mbox{and, then, select the minimum $r$}\cr
& \Quads{1}\mbox{corresponding to this $s$};\cr
\phitm_{w'}(\pair{p, q}), &\mbox{if }(\forall \pair{r,s} <
q)[\phitm_w(\pair{r,s})\converges \not= q];\cr
\diverges, &\mbox{otherwise}.\cr
\end{cases}
\end{equation}

\begin{clm}\label{claim:c-total}
$\PA$ proves $\ll \ \phitm_c$ is total$\ \gg$.
\end{clm}
\Prf{of Claim~\ref{claim:c-total}}
In general proofs of recursion theorems, especially
including Lemma \ref{mixed-rt} above,
are so simple that their proofs can be carried out in $\PA$ (albeit
tediously).
\Qed{Claim~\ref{claim:c-total}}

\begin{clm}\label{claim:w'-mono-increasing}
For all $x > 0$, if $\phitm_{w'}(x-1)\converges \AND
\phitm_{w'}(x)\converges$, then $\phitm_{w'}(x-1) < \phitm_{w'}(x)$.
\end{clm}
\Prf{of Claim~\ref{claim:w'-mono-increasing}}
This follows directly from (\ref{eqn:obfuc-comp-w'}).
\Qed{Claim~\ref{claim:w'-mono-increasing}}

\begin{clm}\label{claim:w'-total}
$\phitm_{w'}$ is total, and $\PA$ proves that.
\end{clm}
\Prf{of Claim~\ref{claim:w'-total}}
Assume for induction that for arbitrarily-fixed $\pair{p, q}$, for all
$x < \pair{p, q}$, we have $\phitm_{w'}(x)\converges$.
By Claim~\ref{claim:c-total}, $\phitm_c(\pair{p, q})\converges$. By
Claim~\ref{claim:pad-works} the output of $\pad$ is always even and
greater
than the input. Thus, repeatedly applying $\pad$ to $\phitm_c(\pair{p,
q})$ will result in $v$ having an even value, $v \not= w$, $v > p$, $v
> q$, and $v > \pair{p, q}$. Furthermore, by the induction assumption,
[if $\pair{p, q} > 0$, then $\phitm_{w'}(\pair{p, q}-1)\converges$],
and padding $\phitm_c(\pair{p, q})$ until $v > \phitm_{w'}(\pair{p,
q}-1)$ is certainly possible. Lastly, by the induction assumption,
checking whether for $(\pair{r, s} = v)$ there does or does not
$(\exists x < \pair{p, q})[\phitm_{w'}(x) = r \Or \phitm_{w'}(x) = s]$
is algorithmically testable.

Thus, by induction, $\phitm_{w'}$ is total. The above inductive proof
is accessible to $\PA$, thus $\PA$ proves it.
\Qed{Claim~\ref{claim:w'-total}}

\begin{clm}\label{claim:range-w-subset-range-w'}
$\rho(\phitm_{w}) \subseteq \rho(\phitm_{w'})$. Furthermore, $\PA$
proves this.
\end{clm}
\Prf{of Claim~\ref{claim:range-w-subset-range-w'}}
For each input $x$, exactly one of the three clauses of
(\ref{eqn:obfuc-comp-w}) must hold. If the first clause holds, then
whatever value $\phitm_w(x)$ has must be a value that was in the range
of $\phitm_w$ on some other input; by implicit application of
Lemma~\ref{lem-pairing}, this recursion must bottom out at
some value, and that value must come from some other clause. 
If the
second clause holds, then $\phitm_w(x) = \phitm_{w'}(x)$, and thus is a
value in the range of $\phitm_{w'}$. Lastly, if the third clause holds,
then $\phitm_w(x)\diverges$, and no value is added to the range of
$\phitm_w$. Thus, all values in the range of $\phitm_w$ have some $x$
such that the second clause of (\ref{eqn:obfuc-comp-w}) holds for
$\phitm_w(x)$, and therefore, the range of $\phitm_w$ is a (potentially
proper) subset of the range of $\phitm_{w'}$.
$\PA$ can handle the preceding argument.
\Qed{Claim~\ref{claim:range-w-subset-range-w'}}

\begin{clm}\label{claim:range-w-checkable-below-q-is-enough}
If $q \in \rho(\phitm_w)$, then
$(\exists \pair{r, s} < q)[\phitm_w(\pair{r, s})\converges = q]$.
Furthermore, $\PA$ proves this.
\end{clm}
\Prf{of Claim~\ref{claim:range-w-checkable-below-q-is-enough}}
Fix arbitrary q, assume that $q \in \rho(\phitm_w)$. Per the
\emph{proof} of Claim~\ref{claim:range-w-subset-range-w'} and the
second clause of (\ref{eqn:obfuc-comp-w}), there then exists $x$ such
that $\phitm_w(x) = q = \phitm_{w'}(x)$. By
Claims~\ref{claim:w'-mono-increasing} and~\ref{claim:w'-total}, 
for $x$ such
that $\phitm_{w'}(x) = q$, it must be the case that $x \leq q$.
By~(\ref{eqn:obfuc-comp-w'}), $x \neq q$.
Thus, $(\exists x < q)[\phitm_w(x) = q]$; as this follows from the
assumption that $q \in \rho(\phitm_w)$, it therefore follows that if $q
\in \rho(\phitm_w)$, $(\exists x < q)[\phitm_w(x) = q]$. 
This proves the claim except for the part about $\PA$.
Lastly, $\PA$ can handle the just prior reasoning.
\Qed{Claim~\ref{claim:range-w-checkable-below-q-is-enough}}

\begin{clm}\label{claim:s_not_in_range_phitm_w}
For all $p$, if $(\exists \pair{r, s} < p)[\phitm_w(\pair{r, s}) = p]$,
then the minimum such $s$ is \emph{not} in the range of $\phitm_w$.
\end{clm}
\Prf{of Claim~\ref{claim:s_not_in_range_phitm_w}}
Assume by way of contradiction otherwise. Fix the least counterexample 
$p$ to the claim, and fix
the minimum $s$ corresponding to that $p$. 
By the assumption by way of contradiction, $(\exists
\pair{r',s'})[\phitm_w(\pair{r', s'}) = s]$. Fix the minimum such $s'$ and
the corresponding $r'$ such that $\pair{r', s'}$ is minimum. By
Claim~\ref{claim:range-w-checkable-below-q-is-enough}, $\pair{r', s'} <
s$.
By the first clause of (\ref{eqn:obfuc-comp-w}), $\phitm_w(\pair{r, s})
= \phitm_w(\pair{\phitm_w(\pair{r, r'}), s'})$. Since $s$ is minimum as
indicated above and $s' < s$ (by Lemma~\ref{lem-pairing}), it
follows that $\pair{\phitm_w(\pair{r, r'}), s'} \geq p$.
\begin{description}
\item[Case one] The first clause of (\ref{eqn:obfuc-comp-w}) holds for input
$\pair{\phitm_w(\pair{r, r'}), s'}$.
Then $s'$ is also a counterexample to 
Claim~\ref{claim:s_not_in_range_phitm_w}, 
which, since $s' < p$, contradicts $p$'s minimality.

\item[Case two] The second clause of (\ref{eqn:obfuc-comp-w}) holds for input
$\pair{\phitm_w(\pair{r, r'}), s'}$. Then
\begin{displaymath}
\phitm_w(\pair{\phitm_w(\pair{r, r'}), s'}) =
\phitm_{w'}(\pair{\phitm_w(\pair{r, r'}), s'}).
\end{displaymath}
Then, by (\ref{eqn:obfuc-comp-w'}), it follows that
\begin{displaymath}
\phitm_w(\pair{\phitm_w(\pair{r, r'}), s'}) > \pair{\phitm_w(\pair{r,
r'}), s'},
\end{displaymath}
and thus, by substitution, it follows that $p > \pair{\phitm_w(\pair{r,
r'}), s'}$; this contradicts
\begin{displaymath}
\pair{\phitm_w(\pair{r, r'}), s'} \geq p.
\end{displaymath}

\item[Case three] The third clause of (\ref{eqn:obfuc-comp-w}) holds for
input $\pair{\phitm_w(\pair{r, r'}), s'}$; this contradicts
$\phitm_w(\pair{r, s})\converges = p$.
\end{description}
\noindent All cases lead to a contradiction; therefore, the claim holds.
\Qed{Claim~\ref{claim:s_not_in_range_phitm_w}}

\begin{clm}\label{claim:w-total}
$\phitm_w$ is total; furthermore, $\PA$ proves this.
\end{clm}
\Prf{of Claim~\ref{claim:w-total}}
Assume for induction that for arbitrarily-fixed $\pair{p, q}$, for all
$x < \pair{p, q}$, $\phitm_w(x)\converges$.
It follows from inequalities about pairing (Lemma~\ref{lem-pairing})
and the induction assumption
that $(\forall \pair{r,s} < q)[\phitm_w(\pair{r,s})\converges]$. Thus,
for input $\pair{p, q}$, the first clause of (\ref{eqn:obfuc-comp-w})
holds if $(\exists \pair{r,s} < q)[q = \phitm_w(\pair{r,s})]$ and the
second clause of (\ref{eqn:obfuc-comp-w}) holds if $(\neg\exists
\pair{r,s} < q)[q = \phitm_w(\pair{r,s})]$; thus, the third clause of
(\ref{eqn:obfuc-comp-w}) does not hold for input $\pair{p, q}$.

If the second clause of (\ref{eqn:obfuc-comp-w}) holds for input
$\pair{p, q}$, then by Claim~\ref{claim:w'-total}, $\phitm_w(\pair{p,
q})\converges$.

If the first clause of (\ref{eqn:obfuc-comp-w}) holds for input
$\pair{p, q}$,
then
\begin{displaymath}
\phitm_w(\pair{p, q})\converges \Iff \phitm_w(\pair{\phitm_w(\pair{p,
r}), s})\converges,
\end{displaymath}
where $\pair{r, s} < q$, $q = \phitm_w({\pair{r, s}})$, $s$ is minimum
such that that is the case, and $r$ is minimum corresponding to $s$. It
follows that both $r < q$ and $s < q$. From the fact that $r < q$, it
follows that $\phitm_w(\pair{p, r})\converges$ by the inductive
assumption. As $s < q$, and $s \not\in \rho(\phitm_w)$ (by
Claim~\ref{claim:s_not_in_range_phitm_w}), it then follows from the
induction hypothesis that $(\forall \pair{r', s'} <
s)[\phitm_w(\pair{r', s'})\converges \not= s]$; thus, the second clause
of (\ref{eqn:obfuc-comp-w}) holds for $\phitm_w(\pair{\phitm_w(\pair{p,
r}), s})$, which is therefore defined by Claim~\ref{claim:w'-total}.

From the inductive assumption, then, $\phitm_w(\pair{p, q})\converges$.
By induction, $\phitm_w$ is total. $\PA$ can prove this as well.
\Qed{Claim~\ref{claim:w-total}}

\begin{clm}\label{claim:psi_w_is_phitm_w}
$\zeta_w = \phitm_w$, and $\PA$ proves this.
\end{clm}
\Prf{of Claim~\ref{claim:psi_w_is_phitm_w}}
The equality follows immediately from the first disjunct of the
disjunction in the first clause of (\ref{eqn:obfuc-comp-psi}); $\PA$
can prove this as well.
\Qed{Claim~\ref{claim:psi_w_is_phitm_w}}

\begin{clm}\label{claim:phitm_phitm_w'_p_q_equals_psi_p_psi_q}
$(\forall p, q)[\phitm_{\phitm_{w'}(\pair{p, q})} = \zeta_p \circ
\zeta_q]$; furthermore, $\PA$ proves this.
\end{clm}
\Prf{of Claim~\ref{claim:phitm_phitm_w'_p_q_equals_psi_p_psi_q}}
The claim follows from (\ref{eqn:obfuc-comp-w'},
\ref{eqn:obfuc-comp-c}) and Claim~\ref{claim:w'-total} above.
\Qed{Claim~\ref{claim:phitm_phitm_w'_p_q_equals_psi_p_psi_q}}

\begin{clm}\label{claim:w_equivalent_w'}
$(\forall p, q)[\phitm_{\phitm_{w'}(\pair{p, q})} =
\phitm_{\phitm_{w}(\pair{p, q})}]$; furthermore, $\PA$ proves this.
\end{clm}
\Prf{of Claim~\ref{claim:w_equivalent_w'}}
Assume for induction that for arbitrary $p, q$, for all $q' < q$, for
all $p'$, $\phitm_{\phitm_{w'}(\pair{p', q'})} =
\phitm_{\phitm_{w}(\pair{p', q'})}$.

Case one: The first clause of (\ref{eqn:obfuc-comp-w}) holds for
$\pair{p, q}$.
 Then, there exists $(\pair{r, s} < q)[q = \phitm_w(\pair{r, s})]$; fix
 $r$ and
 $s$ in accordance with that clause. 
Then fix $r$ and $s$ in accordance with that clause; thus,
  $\pair{r,s} < q = \phi_w(\pair{r,s})$.
It follows that
$\phitm_{w}(\pair{p,
q}) = \phitm_{w}(\pair{\phitm_w(\pair{p, r}), s})$. Since $\pair{r, s}
< q$, it
follows by Lemma~\ref{lem-pairing}
that $r < q$ and $s < q$; thus, by the induction assumption,
$\phitm_{\phitm_{w}(\pair{\phitm_w(\pair{p, r}), s})} =
\phitm_{\phitm_{w'}(\pair{\phitm_w(\pair{p, r}), s})}$ and
$\phitm_{\phitm_w(\pair{p, r})} = \phitm_{\phitm_{w'}(\pair{p, r})}$.
By
Claim~\ref{claim:phitm_phitm_w'_p_q_equals_psi_p_psi_q},
$\phitm_{\phitm_{w'}(\pair{p, r})} = \zeta_p \circ \zeta_r$, and
$\phitm_{\phitm_{w'}(\pair{\phitm_w(\pair{p, r}), s})} =
\zeta_{\phitm_w(\pair{p, r})} \circ \zeta_s$. By
Claim~\ref{claim:range-w-checkable-below-q-is-enough} and the second
disjunct
of the first clause of (\ref{eqn:obfuc-comp-psi}),
$\zeta_{\phitm_w(\pair{p, r})} = \phitm_{\phitm_w(\pair{p, r})}$. By
repeated
substitutions, it follows that $\phitm_{\phitm_{w}(\pair{p, q})} =
\phitm_{\phitm_{w}(\pair{\phitm_w(\pair{p, r}), s})} =
\phitm_{\phitm_{w'}(\pair{\phitm_w(\pair{p, r}), s})} =
\zeta_{\phitm_w(\pair{p, r})} \circ \zeta_s = \phitm_{\phitm_w(\pair{p,
r})}
\circ \zeta_s = \phitm_{\phitm_{w'}(\pair{p, r})} \circ \zeta_s =
\zeta_p \circ
\zeta_r \circ \zeta_s$. It follows from
Claim~\ref{claim:phitm_phitm_w'_p_q_equals_psi_p_psi_q} that
$\phitm_{\phitm_{w'}(\pair{p, q})} = \zeta_p \circ \zeta_q$. As $q =
\phitm_w(\pair{r, s})$ and is thus in $\rho(\phitm_w)$, it follows from
(\ref{eqn:obfuc-comp-psi}) that $\zeta_q = \phitm_q$, and by
substitution it
follows that $\phitm_q = \phitm_{\phitm_w(\pair{r, s})}$. From the
induction
assumption, it follows that $\phitm_{\phitm_w(\pair{r, s})} =
\phitm_{\phitm_{w'}(\pair{r, s})}$, and then from
Claim~\ref{claim:phitm_phitm_w'_p_q_equals_psi_p_psi_q} it follows that
$\phitm_{\phitm_{w'}(\pair{r, s})} = \zeta_r \circ \zeta_s$; and thus
it
follows that $\phitm_{\phitm_{w'}(\pair{p, q})} = \zeta_p \circ \zeta_r
\circ
\zeta_s$. It then immediately follows that
$\phitm_{\phitm_{w'}(\pair{p, q})} = \phitm_{\phitm_{w}(\pair{p,
q})}$.

Case two: The second clause of (\ref{eqn:obfuc-comp-w}) holds for
$\pair{p, q}$. Then, it immediately follows from that clause that
$\phitm_{\phitm_{w'}(\pair{p, q})} = \phitm_{\phitm_{w}(\pair{p,
q})}$.

Case three: The third clause of (\ref{eqn:obfuc-comp-w}) holds for
$\pair{p, q}$. By Claim~\ref{claim:w-total}, this case is impossible.

In all possible cases, $\phitm_{\phitm_{w'}(\pair{p, q})} =
\phitm_{\phitm_{w}(\pair{p, q})}$; therefore, by induction,
\begin{displaymath}
(\forall p, q)[\phitm_{\phitm_{w'}(\pair{p, q})} =
\phitm_{\phitm_{w}(\pair{p, q})}].
\end{displaymath}
$\PA$ can handle the above reasoning, and thus the claim follows.
\Qed{Claim~\ref{claim:w_equivalent_w'}}

\begin{clm}\label{claim:PA-proves-c-witness-composition}
$\PA \proves \ \ll w$ is a witness to composition in $\zeta\gg$.
\end{clm}
\Prf{of Claim~\ref{claim:PA-proves-c-witness-composition}}
This claim follows directly from Claims~\ref{claim:w-total},
\ref{claim:psi_w_is_phitm_w},
\ref{claim:phitm_phitm_w'_p_q_equals_psi_p_psi_q},
and~\ref{claim:w_equivalent_w'}.
\Qed{Claim~\ref{claim:PA-proves-c-witness-composition}}

\begin{clm}\label{claim:range-w-even}
All values in the range of $\phitm_w$ are even.
\end{clm}
\Prf{of Claim~\ref{claim:range-w-even}}
By (\ref{eqn:obfuc-comp-w'}), the range of $\phitm_{w'}$ consists of
only even numbers. The claim follows from that fact and
Claim~\ref{claim:range-w-subset-range-w'}.
\Qed{Claim~\ref{claim:range-w-even}}

In the remainder of this proof, we will need the ability to treat
programs in
the range of $\phitm_w$ --- which are now known to be compositions of
other
programs --- as a sequence of such compositions. Furthermore, we need
the
ability to take potentially-long such sequences and pull off a single
element
from the front or back, and recompose the rest of the sequence to get
another
$\zeta$-program. In order to do this, we introduce $\chain$ and
$\unchain$;
$\chain$ takes the code number of a $\zeta$-program and outputs a
\emph{sequence} of $\zeta$-programs such that if the sequence is
composed in
order, the computed partial function is equivalent to the partial
function computed by
the input $\zeta$-program; while $\unchain$ takes a non-empty sequence
of
$\zeta$-programs and outputs a single $\zeta$-program which is
equivalent to composing the sequence in order.

\begin{equation}\label{eqn:obfuc-comp-chain}
\chain(p) \eqdef \begin{cases}
\chain(r), s, & \mbox{if }
		 (\exists \pair{r, s} < p)[\phitm_w(\pair{r, s}) = p],
		 \cr
&\Quads{1} \mbox{then select the minimum such $s$}\cr
&\Quads{1}\mbox{and, then, select the minimum such $r$}\cr &
\Quads{1}\mbox{corresponding to this selected s;}\cr
p, & \mbox{otherwise.}\cr
\end{cases}
\end{equation}

\noindent When $\mathbf{v}$ is a sequence of more than one element, the terminology
employed just
below,
all-but-last$(\mathbf{v})$ and last$(\mathbf{v})$, is self-explanatory.
\begin{equation}\label{eqn:obfuc-comp-unchain}
\unchain(\mathbf{v}) \eqdef \begin{cases}
\mbox{that element}, & \mbox{if }\mathbf{v}\mbox{ is a sequence}\cr
& \Quads{1}\mbox{of one element;}\cr
\phitm_w(\pair{\unchain(\mbox{all-but-last}(\mathbf{v})), \mbox{last}(\mathbf{v})}), &
\mbox{otherwise.}\cr
\end{cases}
\end{equation}

\begin{clm}\label{claim:psi_chain_x_is_psi_x}
For $p_0, p_1, \ldots, p_n = \chain(x)$, $\zeta_x = \zeta_{p_0} \circ
\zeta_{p_1} \circ \ldots \circ \zeta_{p_n}$.
\end{clm}
\Prf{of Claim~\ref{claim:psi_chain_x_is_psi_x}}
Let $x$ be an arbitrary value, and assume for induction that for all $y
< x$, the claim holds for $\chain(y)$.
\begin{description}
\item[Case one] $\chain(x) = x$. It immediately follows that the claim holds
for $\chain(x)$ in this case.

\item[Case two] $\chain(x) = \chain(p), q$. Then $\phitm_w(\pair{p, q}) = x$,
$x > p$, and  $x > q$.

Then $p$ and $q$ are such that the second clause of (\ref{eqn:obfuc-comp-w})
holds
on input $\pair{p, q}$, per Claim~\ref{claim:s_not_in_range_phitm_w}.
Then $x =
\phitm_{w'}(\pair{p, q})$, 
from which it follows that
$\phitm_x = \zeta_p \circ \zeta_q$. Furthermore, since $x$ is in the
range of
$\phitm_w$, $\zeta_x = \phitm_x$ by the second disjunct of the first
clause of (\ref{eqn:obfuc-comp-psi}). From the just previous reasoning
and the induction assumption, the claim also holds for $\chain(x)$ in
this case.
\end{description}
\noindent By the induction, the claim holds for all $x$.
\Qed{Claim~\ref{claim:psi_chain_x_is_psi_x}}

\begin{clm}\label{claim:elems-of-chain-are-not-range-w}
For any $p$, none of the elements of $\chain(p)$ are in
$\rho(\phitm_w)$.
\end{clm}
\Prf{of Claim~\ref{claim:elems-of-chain-are-not-range-w}}
Assume by induction that for all $p' < p$, the elements of $\chain(p')$
are each not in $\rho(\phitm_w)$.
\begin{description}
\item[Case one] The first clause of (\ref{eqn:obfuc-comp-chain}) holds for
$p$. Let $r, s$ be per that clause; by Lemma~\ref{lem-pairing}, it
follows that $r < p$. Therefore, by the induction assumption, the
elements of $\chain(r)$ are each not in $\rho(\phitm_w)$. By
Claim~\ref{claim:s_not_in_range_phitm_w}, $s$ is also not in
$\rho(\phitm_w)$; thus, each element of $\chain(p)$ is not in
$\rho(\phitm_w)$.

\item[Case two] The second clause of (\ref{eqn:obfuc-comp-chain}) holds for
$p$. By Claim~\ref{claim:range-w-checkable-below-q-is-enough}, if $p$
is in the range of $\phitm_w$, then $(\exists\pair{r, s} <
p)[\phitm_w(\pair{r, s})\converges = p]$; thus, $p \not\in
\rho(\phitm_w)$, and the only element of $\chain(p)$ is not in
$\rho(\phitm_w)$.
\end{description}
\noindent By induction, the claim holds for all $p$.
\Qed{Claim~\ref{claim:elems-of-chain-are-not-range-w}}

\begin{clm}\label{claim:unchain_chain_x_is_x}
For all $x$, $\unchain(\chain(x)) = x$.
\end{clm}
\Prf{of Claim~\ref{claim:unchain_chain_x_is_x}}
If $\chain(x)$ is a sequence of one element, $\unchain(\chain(x))$ is
just $x$ per the first clause of (\ref{eqn:obfuc-comp-unchain}).
Let $n$ be some number greater than $1$. Assume for induction that for
all $x$ such that $\chain(x)$ is a sequence of length no more than
$n-1$, $\unchain(\chain(x))=x$. Then let $x$ be such that $\chain(x)$
is a sequence of length n. That sequence, for some $r, s$, dependent on
$x$, is $\chain(r), s$. Furthermore, by the first clause of
(\ref{eqn:obfuc-comp-chain}), $x = \phitm_w(\pair{r, s})$. In such a
sequence, for $\unchain(\chain(r), s)$, the second clause of
(\ref{eqn:obfuc-comp-unchain}) holds, and outputs the result of
$\phitm_w(\pair{\unchain(\chain(r)),s})$ --- by inductive assumption,
it follows that $\unchain(\chain(s))$ is s, from which it follows that
the output of $\unchain(\chain(x))$ is $\phitm_w(\pair{r, s})$ --- which
is already known to be $x$.
Thus, the claim follows by induction on the length of $\chain(x)$.
\Qed{Claim~\ref{claim:unchain_chain_x_is_x}}

\begin{clm}\label{claim:partial-unchain-is-ok}
For any sequence $p_0,p_1,..,p_n$ of length at least two such that
there exists
$p$ such that $\chain(p) = p_0,p_1,\ldots,p_n$, $\zeta_p = \zeta_{p_0}
\circ \zeta_{\unchain(p_1,\ldots,p_n)}$.
\end{clm}
\Prf{of Claim~\ref{claim:partial-unchain-is-ok}}
Let $p' = \unchain(p_1,..,p_n)$. By (\ref{eqn:obfuc-comp-chain},
\ref{eqn:obfuc-comp-unchain}) above and by
Claim~\ref{claim:elems-of-chain-are-not-range-w}, it follows that
$\chain(p') =
p_1,..,p_n$. Therefore, by Claim~\ref{claim:psi_chain_x_is_psi_x},
$\zeta_{p'}
= \zeta_{p_1} \circ \ldots \circ \zeta_{p_n}$. Thus, $\zeta_{p_0} \circ
\zeta_{\unchain(p_1,..,p_n)} = \zeta_{p_0} \circ \zeta_{p_1} \circ
\ldots \circ
\zeta_{p_n}$; by Claim~\ref{claim:psi_chain_x_is_psi_x}, this is
exactly $\zeta_p$.
\Qed{Claim~\ref{claim:partial-unchain-is-ok}}

\begin{clm}\label{claim:w-increases-many-inputs}
For all $x$ such that $[x$ is odd \emph{or} $x$ in $\rho(\phitm_w)]$,
$\zeta_w(x) > x$, and, if $x$ is in $\rho(\phitm_w)$, then there does
not exist
$y \not= x$ such that $\zeta_w(y) = \zeta_w(x)$.
\end{clm}
\Prf{of Claim~\ref{claim:w-increases-many-inputs}}
Suppose $x$ is odd. For $\pair{p, q} = x$, $q$ must be odd --- per 
Lemma~\ref{lem-pairing}.
As such, because
$\phitm_w$ only outputs even values, $q$ cannot be in the range of
$\phitm_w$,
and thus $\phitm_w(x) = \phitm_{w'}(x)$ by (\ref{eqn:obfuc-comp-w}). By
(\ref{eqn:obfuc-comp-w'}), $\phitm_{w'}(x) > x$. By
Claim~\ref{claim:psi_w_is_phitm_w}, $\zeta_w = \phitm_w$, and thus
$\zeta_w(x) > x$.

Suppose $x \in \rho(\phitm_w)$. For $\pair{p, q} = x$, by 
Lemma~\ref{lem-pairing}, $p \leq x$ and $q \leq x$. By
(\ref{eqn:obfuc-comp-w'}) and Claim~\ref{claim:w'-total}, there does not exist $y < x$ such that
$\phitm_{w'}(y) =$ either $p$ or $q$. 
By~(\ref{eqn:obfuc-comp-w'}) and Claim~\ref{claim:w'-total}, $\phi_{w'}(0) > 0$;
from this and
Claims~\ref{claim:w'-mono-increasing} and~\ref{claim:w'-total}, it follows that
there does not exist $y$ such that $\phitm_{w'}(y) = $ either $p$ or
$q$. 
Therefore, by Claim~\ref{claim:range-w-subset-range-w'}, $q
\not\in \rho(\phitm_{w'})$, from which it follows, by
(\ref{eqn:obfuc-comp-w}), $\phitm_w(x) = \phitm_{w'}(x)$. As $p \not\in
\rho(\phitm_{w'})$ and thus $p \not\in \rho(\phitm_w)$, it follows by 
pairing being 1-1 that
there does not exist any $p', r, s$ such that
$\pair{\phitm_w(\pair{p',r}), s} = x$. 

Therefore, there exists no $y$
such that both the first clause of (\ref{eqn:obfuc-comp-w}) holds on
input $y$ and $\phitm_w(y) = \phitm_w(x)$. By Claims
\ref{claim:w'-mono-increasing} and \ref{claim:w'-total}, there exists no
$y \not= x$ such that $\phitm_{w'}(x) = \phitm_{w'}(y)$; from this
and the fact that $\phi_w(x) = \phi_{w'}(x)$,
it follows that for $y \not= x$, if the second clause of
(\ref{eqn:obfuc-comp-w}) holds on input $y$, $\phitm_w(y) \not= \phitm_w(x)$.
Therefore, there does not exist $y \not= x$
such that $[\phitm_w(x) = \phitm_w(y)]$.
Furthermore, from (\ref{eqn:obfuc-comp-w'}) and the fact that $\phitm_w(x) =
\phitm_{w'}(x)$, it follows that $\phitm_w(x) > x$.

Thus, the claim follows.
\Qed{Claim~\ref{claim:w-increases-many-inputs}}

\begin{clm}\label{claim:T-proves-equality-implies-p-increasing}
If, for some $x_0$, $\T \proves_{x_0} \ \ll (\exists r, t \not= r)[\zeta_r =
\zeta_t]\gg$, then
$(\forall p)(\forall x \geq x_0 \mid x $ is odd $ \Or x \in
\rho(\phitm_w))[\zeta_p(x)\converges \AND [\zeta_p(x)$ odd $\Or
\zeta_p(x) \in
\rho(\phitm_w)] \AND \zeta_p(x) > x]$.
\end{clm}
\Prf{of Claim~\ref{claim:T-proves-equality-implies-p-increasing}}
Assume by way of contradiction that the claim does not hold. Fix least
$p$ and
least corresponding $x \geq x_0$ such that $[x$ is odd $ \Or x \in
\rho(\phitm_w)]$ so that either $\zeta_p(x)\diverges$, or $\zeta_p(x)$
even and
not in the range of $\phitm_w$, or $\zeta_p(x) \leq x$.
\begin{description}
\item[Case one] $p$ is $w$. Then, by Claims~\ref{claim:w-total},
\ref{claim:w-increases-many-inputs} and~\ref{claim:psi_w_is_phitm_w}, a
contradiction follows immediately.

\item[Case two] $p$ is not $w$ and not in $\rho(\phitm_w)$. Then, by
(\ref{eqn:obfuc-comp-psi}), $\zeta_p(x) = \myprime(p+1)^x$, which is an
odd value greater than $x$. A contradiction follows immediately.

\item[Case three] $p$ is in $\rho(\phitm_w)$. Let $p_0, \ldots, p_n$ be
$\chain(p)$ and
let $p' = \unchain(p_0, \ldots, p_{n-1})$.\footnote{
By (\ref{eqn:obfuc-comp-chain}), for all $x$ in the range of $\phi_w$,
the length of $\chain(x) \geq 2$.}
By Claim~\ref{claim:elems-of-chain-are-not-range-w}, $p_n \not\in
\rho(\phitm_w)$.
By Claim~\ref{claim:unchain_chain_x_is_x} and
(\ref{eqn:obfuc-comp-chain}),
$\phitm_w(\pair{p', p_n}) = p$. 
From (\ref{eqn:obfuc-comp-chain}) we have that
$p', p_n, \pair{p', p_n} < p$. Then
by Claims~\ref{claim:phitm_phitm_w'_p_q_equals_psi_p_psi_q} 
and~\ref{claim:w_equivalent_w'} and (\ref{eqn:obfuc-comp-psi}), it follows
that $\zeta_p =
\zeta_{p'} \circ \zeta_{p_n}$. By the second clause of
(\ref{eqn:obfuc-comp-w}), $\phitm_w(\pair{p', p_n}) =
\phitm_{w'}(\pair{p', p_n})$, thus $p > p'$ by
(\ref{eqn:obfuc-comp-w'}).

\begin{quotation}
Subcase one: $p_n$ is not $w$. Then, $\zeta_{p_n}(x) =
\myprime(p_n+1)^x$, which
is both odd and greater than $x$; thus, by the assumption that $p$ is
the least
value such that the claim does not hold,
$\zeta_{p'}(\myprime(p_n+1)^x)$ is
greater than $\myprime(p_n+1)^x$, which is greater than $x$.
$\zeta_{p'}(\myprime(p_n+1)^x)$ is also either odd or in
$\rho(\phitm_w)$; as
$\zeta_p(x)$ is that value, contradiction follows immediately.
\end{quotation}

\begin{quotation}
Subcase two: $p_n$ is $w$. 
Then, by (\ref{eqn:obfuc-comp-psi}),
$\zeta_{p_n}(x)$ is in the range of $\phitm_w$, moreover, by
Claim~\ref{claim:w-increases-many-inputs}, it is $> x$.
By the assumption that $p$ is the
least value such that the claim does not hold,
$\zeta_{p'}(\zeta_{p_n}(x))$ is
greater than $\zeta_{p_n}(x)$, which, as previously shown, is greater
than $x$.
$\zeta_{p'}(\zeta_{p_n}(x))$ is also either odd or in $\rho(\phitm_w)$;
as that
value is equal to $\zeta_p(x)$ --- a contradiction follows immediately.
\end{quotation}
\end{description}
In all cases, a contradiction follows; thus, the claim holds.
\Qed{Claim~\ref{claim:T-proves-equality-implies-p-increasing}}

\begin{clm}\label{claim:T-proves-equality-implies-infinitely-not-equal}
If, for some $x_0$, $\T \proves_{x_0} \ \ll (\exists r, t \not= r)[\zeta_r =
\zeta_t]\gg$, then
$(\forall r', t' \not= r')(\existsio x)[\zeta_{r'}(x) \not=
\zeta_{t'}(x)]$.
\end{clm}
\Prf{of
Claim~\ref{claim:T-proves-equality-implies-infinitely-not-equal}}
Assume by way of contradiction otherwise. Then by the assumption that
$\T$ does
not prove false things, $(\exists p', q' \not= p')[\zeta_{p'} =
\zeta_{q'}]$.
Thus, $(\exists p', q' \not= p')(\forall x \geq 0)[\zeta_{p'}(x) =
\zeta_{q'}(x)]$. Let $p$ be the least number
such that $\chain(p)$ is of minimum
length so
that there are $q$ and $x''$ with $q \not= p$ and $x''$ 
so that $(\forall x \geq
x'')[\zeta_p(x) = \zeta_q(x)]$; then let $q$ be the least corresponding
$q$.
Let $x_0$ be the least number such that $\T \proves_{x_0} \ \ll
(\exists r, t
\not= r)[\zeta_r = \zeta_t]\gg$. Let $x'$ be the least value such that
each of the following hold: $x'$ is in the range of $\phitm_w$, $x'
\geq x_0$, and $x' \geq x''$. By
Claim~\ref{claim:w-increases-many-inputs}, there is a value in the
range of $\phitm_w$ greater than any fixed odd number, and thus such an
$x'$ must exist.

\begin{description}
\item[Case 1] $p = w$.

\begin{quotation}
Subcase 1.1: $q = w$. This immediately contradicts our
assumption that $p \not= q$.
\end{quotation}

\begin{quotation}
Subcase 1.2: $q \not= w$ and $q \not\in \rho(\phitm_w)$. Then,
by
(\ref{eqn:obfuc-comp-psi}), $\zeta_q(x') = \myprime(q+1)^{x'}$, which
is an odd
value, and by  Claim~\ref{claim:range-w-even}, $\zeta_w(x')$ is even; a
contradiction follows immediately.
\end{quotation}

\begin{quotation}
Subcase 1.3: $q \in \rho(\phitm_w)$. Let $q_0, \ldots, q_n =
\chain(q)$. By
Claim~\ref{claim:elems-of-chain-are-not-range-w}, $q_0$ is not in
$\rho(\phitm_w)$. Let $q' = \unchain(q_1, \ldots, q_n)$. By
Claim~\ref{claim:partial-unchain-is-ok}, $\zeta_q = \zeta_{q_0} \circ
\zeta_{q'}$. By
Claim~\ref{claim:T-proves-equality-implies-p-increasing},
$\zeta_{q'}(x') > x'$ and $\zeta_{q'}(x')$ is either odd or in the
range of $\phitm_w$.
\end{quotation}

   \begin{quotation}
   \begin{quotation}
Sub-subcase 1.3.1: $q_0 \not= w$. Then by
(\ref{eqn:obfuc-comp-psi}),
$\zeta_q(x) = \zeta_{q_0}(\zeta_{q'}(x')) =
\myprime(q_0+1)^{\zeta_{q'}(x')}$,
which is an odd value, while, by Claim~\ref{claim:range-w-even},
$\zeta_w(x')$
is even; therefore, $\zeta_q(x') \not= \zeta_p(x')$; a contradiction.
  \end{quotation}
  \end{quotation}

   \begin{quotation}
   \begin{quotation}
Sub-subcase 1.3.2: $q_0 = w$. Therefore, since
$\zeta_{q_0}(\zeta_{q'}(x')) = \zeta_q(x') = \zeta_p(x')$,
$\zeta_{q_0}(\zeta_{q'}(x')) = \zeta_w(\zeta_{q'}(x'))$, and
$\zeta_p(x') =
\zeta_w(x')$, it follows that $\zeta_q(x') = \zeta_w(x')$. However,
$\zeta_{q'}(x') > x'$, and by
Claim~\ref{claim:w-increases-many-inputs}, it
follows from the fact that $\zeta_w(x') = \zeta_w(\zeta_{q'}(x'))$ that
$x' =
\zeta_{q'}(x')$ --- a contradiction.
  \end{quotation}
  \end{quotation}

\item[Case 2] $p \not= w$ and $p \not\in \rho(\phitm_w)$.

\begin{quotation}
Subcase 2.1: $q = w$. The same argument holds as for 
Case~1.2,
interchanging $p$ and $q$.
\end{quotation}

\begin{quotation}
Subcase 2.2: $q \not= w$ and $q \not\in \rho(\phitm_w)$. Then
by
(\ref{eqn:obfuc-comp-psi}), $\zeta_p(x') = \myprime(p+1)^{x'}$ and
$\zeta_q(x') = \myprime(q+1)^{x'}$, thus $\myprime(p+1)^{x'} =
\myprime(q+1)^{x'}$, which implies that $p = q$, a contradiction.
\end{quotation}

\begin{quotation}
Subcase 2.3: $q \in \rho(\phitm_w)$. Let $q_0, \ldots, q_n =
\chain(q)$;
then let $q' = \unchain(q_1, \ldots, q_n)$. By
Claim~\ref{claim:partial-unchain-is-ok}, $\zeta_q = \zeta_{q_0} \circ
\zeta_{q'}$.  By
Claim~\ref{claim:T-proves-equality-implies-p-increasing}, for
all $x \geq x'$, $\zeta_{p'}(x)$ is greater than $x'$ and either odd or
in the range of $\phitm_w$.
By Claim~\ref{claim:elems-of-chain-are-not-range-w},
$q_0 \notin \rho(\phi^{TM}_w)$.
\end{quotation}

\begin{quotation}
\begin{quotation}
Sub-subcase 2.3.1: $q_0 = w$. Then $\zeta_q(x') =
\zeta_w(\zeta_{q'}(x'))$, which is even per
Claim~\ref{claim:range-w-even}. By
the second clause of (\ref{eqn:obfuc-comp-psi}), $\zeta_p(x')$ is odd;
a
contradiction to $\zeta_p(x') = \zeta_q(x')$ follows.
\end{quotation}
\end{quotation}

\begin{quotation}
\begin{quotation}
Sub-subcase 2.3.2: $q_0 = p$. Then, as $\zeta_q(x') =
\zeta_{q_0}(\zeta_{q'}(x')$, it follows that $\zeta_q(x') =
\myprime(p+1)^{\zeta_{q'}(x')}$. Likewise, $\zeta_p(x') =
\myprime(p+1)^{x'}$.
Therefore, $\zeta_{q'}(x') = x'$, which is a contradiction to
Claim~\ref{claim:T-proves-equality-implies-p-increasing}.
\end{quotation}
\end{quotation}

\begin{quotation}
\begin{quotation}
Sub-subcase 2.3.3: $q_0 \not= w$ and $q_0 \not= p$. Then, by
(\ref{eqn:obfuc-comp-psi}) and the fact that $\zeta_q(x') =
\zeta_{q_0}(\zeta_{q'}(x'))$, $\zeta_q(x')$ is a power of
$\myprime(q_0+1)$ and
$\zeta_p(x')$ is $\myprime(p+1)^{x'}$; thus $\zeta_q(x') \not=
\zeta_p(x')$, a contradiction.
\end{quotation}
\end{quotation}

\item[Case 3] $p \in \rho(\phitm_w)$.

\begin{quotation}
Subcase 3.1: $q = w$. The same argument holds as for Case~1.3,
interchanging $p$ and $q$.
\end{quotation}

\begin{quotation}
Subcase 3.2: $q \not= w$ and $q \not\in \rho(\phitm_w)$. The
same argument holds as for case two subcase three, interchanging $p$
and $q$.
\end{quotation}

\begin{quotation}
Subcase 3.3: $q \in \rho(\phitm_w)$. Let $p_0, \ldots, p_m =
\chain(p)$,
let $q_0, \ldots, q_n = \chain(q)$. Let $p' = \unchain(p_1, \ldots, p_m)$ and
$q' =
\unchain(q_1, \ldots, q_n)$. By Claim~\ref{claim:partial-unchain-is-ok},
$\zeta_q =
\zeta_{q_0} \circ \zeta_{q'}$, and $\zeta_p = \zeta_{p_0} \circ
\zeta_{p'}$. By
Claim~\ref{claim:T-proves-equality-implies-p-increasing}, for all $x
\geq x'$,
$\zeta_{p'}(x)$ and $\zeta_{q'}(x)$ are both greater than $x'$ and
either odd or in the range of $\phitm_w$.
\end{quotation}

\begin{quotation}
\begin{quotation}
Sub-subcase 3.3.1: $p_0 = q_0 \not= w$. Then, by
(\ref{eqn:obfuc-comp-psi}), for each $x \geq x'$, $\zeta_{p_0}(x) =
\zeta_{q_0}(x) = \myprime(p_0+1)^x$. Thus, for each $x \geq x'$, for
each $y \geq
x' \AND y \not= x$, $\zeta_{p_0}(x) \not= \zeta_{q_0}(y)$. As it is the
case
that for all $x \geq x'$, $\zeta_{p'}(x)$ and $\zeta_{q'}(x)$ are both
greater
than $x'$, it then follows from the previous equalities and
inequalities about
$p_0$ and $q_0$, as well as the facts that $\zeta_q = \zeta_{q_0} \circ
\zeta_{q'}$, and $\zeta_p = \zeta_{p_0} \circ \zeta_{p'}$, that
$(\forall x
\geq x')[\zeta_{p'}(x) = \zeta_{q'}(x)]$. But since $\chain(p')$ is
shorter than
$\chain(p)$, this contradicts our assumption that $p$ had the shortest
$\chain$
such that there exists $q \not= p$ and $x''$ so that $(\forall x \geq
x'')[\zeta_p(x) = \zeta_q(x)]$.
\end{quotation}
\end{quotation}

\begin{quotation}
\begin{quotation}
Sub-subcase 3.3.2: $p_0 = q_0 = w$. By
Claim~\ref{claim:w-increases-many-inputs} and the fact that for all $x
\geq
x'$, $\zeta_{p'}(x)$ and $\zeta_{q'}(x)$ are both greater than $x'$ and
either
odd or in the range of $\phitm_w$, it follows from reasoning similar to
that of
sub-subcase 3.3.1, that for all $x \geq x'$, $\zeta_{p'}(x) =
\zeta_{q'}(x)$,
which again contradicts our assumption that $p$ had the shortest
$\chain$ such
that there exists $q \not= p$ and $x''$ so that $(\forall x \geq
x'')[\zeta_p(x) = \zeta_q(x)]$.
\end{quotation}
\end{quotation}

\begin{quotation}
\begin{quotation}
Sub-subcase 3.3.3: $p_0 \not= q_0$, and neither $p_0$ nor
$q_0$ is
$w$. Then, by (\ref{eqn:obfuc-comp-psi}) and the facts that $\zeta_p =
\zeta_{p_0} \circ \zeta_{p'}$ and $\zeta_q = \zeta_{q_0} \circ
\zeta_{q'}$, it
follows that $\zeta_p(x')$ is a power of $\myprime(p_0+1)$, and
$\zeta_q(x')$ is a
power of $\myprime(q_0+1)$; thus, $\zeta_p(x') \not= \zeta_q(x')$, a
contradiction.
\end{quotation}
\end{quotation}

\begin{quotation}
\begin{quotation}
Sub-subcase 3.3.4: $p_0 \not= q_0$, and $p_0 = w$. Then, by
Claim~\ref{claim:range-w-even} and the fact that $\zeta_p = \zeta_{p_0}
\circ
\zeta_{p'}$, it follows that $\zeta_p(x')$ is even. By
Claim~\ref{claim:elems-of-chain-are-not-range-w}, $q_0 \not\in
\rho(\phitm_w)$
and $q_0 \not= w$; thus by (\ref{eqn:obfuc-comp-psi}), the fact that
$\zeta_q =
\zeta_{q_0} \circ \zeta_{q'}$, and the fact that $\zeta_{q'}(x') \geq
x'$, it
follows that $\zeta_q(x')$ is odd; a contradiction to the assumption
that
$\zeta_p(x') = \zeta_q(x')$.
\end{quotation}
\end{quotation}

\begin{quotation}
\begin{quotation}
Sub-subcase 3.3.5: $p_0 \not= q_0$, and $q_0 = w$. The same
reasoning holds as for sub-subcase 3.3.4, with p and q interchanged.
\end{quotation}
\end{quotation}
\end{description}
Every case leads to a contradiction; thus the assumption is false, and
the claim follows.
\Qed{Claim~\ref{claim:T-proves-equality-implies-infinitely-not-equal}}

\begin{clm}\label{claim:T-not-proves-equality}
$\T \not\proves \ \ll (\exists r, t \not= r)[\zeta_r = \zeta_t]\gg$.
\end{clm}
\Prf{of Claim~\ref{claim:T-not-proves-equality}}
By Claim~\ref{claim:T-proves-equality-implies-infinitely-not-equal}, if
$\T
\proves \ \ll (\exists r, t \not= r)[\zeta_r = \zeta_t]\gg$, then
$(\forall r',
t' \not= r')(\existsio x)[\zeta_{r'}(x) \not= \zeta_{t'}(x)]$. Thus,
there does
not exist $r, t \not = r$ such that $\zeta_r = \zeta_t$. Since $\T$
does not prove false things, the claim follows.
\Qed{Claim~\ref{claim:T-not-proves-equality}}

\begin{clm}\label{claim:is_acceptable}
$\zeta$ is acceptable.
\end{clm}
\Prf{of Claim~\ref{claim:is_acceptable}}
By Claim~\ref{claim:T-not-proves-equality}, the third disjunct of the
first
clause of (\ref{eqn:obfuc-comp-psi}) is true for all $x$. Thus, for all
$p$ and
$x$, $\zeta_p(x) = \phitm_p(x)$; that is, $\zeta = \phitm$, which is
known to be acceptable.
\Qed{Claim~\ref{claim:is_acceptable}}

The theorem follows immediately from
Claims~\ref{claim:PA-proves-c-witness-composition},
\ref{claim:T-not-proves-equality}, and~\ref{claim:is_acceptable}
\Qed{Theorem~\ref{thm:t-prove-composition-not-equality}}

\section*{Acknowledgement}
We are very grateful for the anonymous
referees' hard work, helpful suggestions, and corrections.


\vspace{-40 pt}
\end{document}